\documentclass[10pt]{article}
 \usepackage{amsmath,amssymb,amsthm,amsfonts}
 \usepackage{epstopdf}
 \RequirePackage[dvips]{graphicx}
 \usepackage[dvips]{epsfig}
 \usepackage{color}
\usepackage[usenames,dvipsnames]{pstricks}

 \def\draw #1 by #2 (#3){
  \vbox to #2{
    \hrule width #1 height 0pt depth 0pt
    \vfill
    \special{picture #3} 
    }
  }

 \def\scaleddraw #1 by #2 (#3 scaled #4){{
  \dimen0=#1 \dimen1=#2
  \divide\dimen0 by 1000 \multiply\dimen0 by #4=
  \divide\dimen1 by 1000 \multiply\dimen1 by #4
  \draw \dimen0 by \dimen1 (#3 scaled #4)}
  }

\newtheorem{theorem}{Theorem}[section]
\newtheorem{example}[theorem]{Example}
\newtheorem{problem}[theorem]{Problem}
\newtheorem{conjecture}[theorem]{Conjecture}
\newtheorem{defin}[theorem]{Definition}
\newtheorem{lemma}[theorem]{Lemma}
\newtheorem{corollary}[theorem]{Corollary}
\newtheorem{remark}[theorem]{Remark}
\newtheorem{nt}{Note}

\newenvironment{pf}{\medskip\noindent{Proof:  \hspace*{-.4cm}}
       \enspace}{\hfill \qed \newline \medskip}

 \newtheorem{proposition}[theorem]{Proposition}
 \newtheorem{rule-def}[theorem]{Rule}

\begin{document}
 \newcommand{\la}{\lambda}
 \newcommand{\si}{\sigma}
 \newcommand{\ol}{1-\lambda}
 \newcommand{\be}{\begin{equation}}
 \newcommand{\ee}{\end{equation}}
 \newcommand{\bea}{\begin{eqnarray}}
 \newcommand{\eea}{\end{eqnarray}}

\author{Gi-Sang Cheon\thanks{Applied Algebra and Optimization Research Center, Department of Mathematics, Sungkyunkwan University, Suwon 16419, Republic of Korea},
 Ji-Hwan Jung\footnotemark[2], Sergey Kitaev\thanks{Department of Computer and Information Sciences, University of Strathclyde, 26 Richmond Street, Glasgow,
  G1 1XH, United Kingdom}\  \ and Seyed Ahmad Mojallal\footnotemark[2]\\
{\footnotesize gscheon@skku.edu},\,{\footnotesize
jh56k@skku.edu},\,{\footnotesize
sergey.kitaev@cis.strath.ac.uk},\,{\footnotesize mojallal@skku.edu}}

\title{{Riordan graphs II: Spectral properties}
\thanks{This work was supported by the National Research Foundation of Korea(NRF) grant funded by the Korea government(MSIP) (2016R1A5A1008055) and the Ministry of Education
(NRF-2016R1A6A3A11930452).}
\date{}}

\maketitle

 \begin{abstract}
 The authors of this paper have used the theory of Riordan matrices to introduce the notion of a Riordan graph in \cite{CJKM}. Riordan graphs are proved to have a number of interesting (fractal) properties, and they are a far-reaching generalization of the well known and well studied Pascal graphs and Toeplitz graphs, and also some other families of graphs. The main focus in \cite{CJKM} is the study of structural properties of families of Riordan graphs obtained from certain infinite Riordan graphs.

In this paper, we use a number of results in~\cite{CJKM} to study
spectral properties of  Riordan graphs. Our studies include, but
are not limited to the spectral graph invariants for Riordan graphs
such as the adjacency eigenvalues, (signless) Laplacian eigenvalues,
nullity, positive and negative inertias, and rank.  We also study
determinants of Riordan graphs, in particular, giving results about determinants of Catalan graphs.

 \bigskip

 \end{abstract}

 \section{Introduction}

The authors of this paper have used the theory of {\em Riordan matrices}  to introduce the notion of a {\em Riordan graph} in  \cite{CJKM}. Riordan graphs are a far-reaching generalization of the well known and well studied
 {\em Pascal graphs} \cite{DQ} and {\em Toeplitz graphs} \cite{Ghorban}, and also some other families of graphs. The Pascal graphs are constructed using {\em Pascal's triangle} modulo~2, and Pascal's triangle
 itself has motivated the appearance of the area of Riordan matrices~\cite{shap}, an active area of research these days. A Toeplitz graph $G=(V, E)$ is a graph with $V =\{1,\ldots,n\}$ and
 $E=\{ij\ |\ |i-j|\in\{t_1,\ldots,t_k\},1\leq t_1<\cdots< t_k\leq n-1\}$.

Riordan graphs are proved to have a number of interesting (fractal) properties \cite{CJKM}, which can be useful, e.g.\ in creating computer networks \cite{DQ} with
certain desirable features, such as
\begin{itemize}
\item the design is to be simple and recursive;
\item there must be a universal vertex adjacent to all others;
\item there must exist several paths between each pair of vertices.
\end{itemize}
Also, Riordan graphs can be useful when designing algorithms to compute values of graph invariants.

The main focus in \cite{CJKM} is the study of {\em structural properties} of families of Riordan graphs obtained from infinite Riordan graphs, which includes  the fundamental {\em Riordan Graph Decomposition Theorem} (see Theorem~\ref{e:th} below) and the generalization of a number of known results for the Pascal graphs.

In this paper, we study {\em spectral properties} of Riordan graphs.
The {\em spectral graph theory} studies relations between graph
properties and the {\em spectrum} (i.e.\ the set of eigenvalues) of
the adjacency matrix, or (signless) Laplacian matrix,  which can be
useful in various contexts. For example,
\begin{itemize}
\item the {\em second largest} eigenvalue of a graph gives information about {\em expansion} and {\em randomness} properties of a graph;
\item the {\em smallest} eigenvalue gives information about the {\em independence} and {\em chromatic numbers} of a graph;
\item  {\em interlacing} of eigenvalues (see the definition in Lemma~\ref{lem31}) gives information about {\em graph substructures};
\item the fact that eigenvalue multiplicities must be integral provides strong restrictions on graph parameters, e.g. in the case of {\em strongly-regular graphs}.
\end{itemize}

In relation to our paper, in \cite{CKM}, several spectral properties
of Pascal graphs were studied by exploring their spectral graph
invariants such as the algebraic connectivity, the first three
largest Laplacian eigenvalues and the nullity. In this paper, we
obtain results on several other spectral graph invariants for Pascal
graphs. More importantly, we develop the spectral theory for many
other classes of Riordan graphs. Our studies include, but are not
limited to spectral graph invariants for Riordan graphs such as the
{\em adjacency eigenvalues}, ({\em signless}) {\em Laplacian
eigenvalues} (see Section~\ref{eigenvalues-sec}), positive and
negative {\em inertias} (see Section~\ref{inertia-sec}), nullity
(see Section~\ref{nullity-sec}), and {\em rank} (see
Theorem~\ref{lem3}). Our key achievement is to give relations for
Laplacian spectral radius and algebraic connectivity of Riordan
graphs in term of given generating functions. Also,  we give a new
version of Rayleigh-Ritz quotient for Riordan graphs.

One of the essential difficulties we faced was the fact that the number of edges in Riordan graphs is not known in the general case. This parameter is normally known when dealing with spectral graph invariants. Still, we were able to obtain non-trivial interesting spectral results for various classes of Riordan graphs. In some cases, we use certain structural results in \cite{CJKM} to achieve our goals.

The paper is organized as follows. In Section~\ref{sec-prelim}, we give a list of necessary notions, notations and known results.
 In Section~\ref{eigenvalues-sec}, we study the eigenvalues and (signless) Laplacian eigenvalues of Riordan graphs.  In
 Section~\ref{vertexdegrees-sec}, vertex degrees in Riordan graphs
 are studied to obtain some bounds for the largest and second
 smallest Laplacian eigenvalues of Riordan graphs.  The {\em Rayleigh-Ritz quotient}
 for Riordan graphs is also studied in Section~\ref{vertexdegrees-sec}.   In Section~\ref{PasCat}, we
 give results on eigenvalues and Laplacian eigenvalues of
 Pascal and Catalan graphs.
In Section~\ref{inertia-sec}, results pertaining to the
positive and negative inertias of a Riordan graph, and their
complements, are given. In Section~\ref{nullity-sec}, we present
several results on the nullity and rank of Riordan graphs. Finally,
in Section~\ref{det-Catalan-sec}, we study determinants of graphs in
certain subclasses of Riordan graphs. In particular, we give results
on determinants of Catalan graphs.

\section{Preliminaries}\label{sec-prelim}

Graphs in this paper are normally on the vertex set
$[n]:=\{1,2,\ldots,n\}$.  Also, we let  $V_o:=\{j\in [n]\ |\, j ~
\mbox{is odd}\}$ and $V_e:=\{j\in [n]\ |\, j ~ \mbox{is even}\}$.
Suppose $G$ is a graph with a vertex set $V$. For $W \subseteq V$,
we denote the subgraph in $G$ induced by $W$ by $\langle W
\rangle_G$, or simply by $\langle W \rangle$ if $G$ is understood
from the context. The neighborhood of a vertex $v$, denoted
by $N_G(v)$, is the set of vertices in a graph $G$ adjacent to $v$.

In this paper, the graphs $N_n$, $K_n$, $K_{m, n-m}\,(1\leq m \leq \lfloor \frac{n}{2}\rfloor)$, and $P_n$ are, respectively, the {\em null graph}, the {\em complete graph},
 the {\em complete bipartite graph} with parts of sizes $m$ and $n-m$, and the {\em path graph} on $n$ vertices. Also, $A\equiv B$ means that matrices $A$ and $B$ equal modulo 2. Similarly, $A(i,j)\equiv B(i,j)$ means the $(i,j)$th element in $A$ equals the $(i,j)$th element in $B$ modulo~2.

\subsection{Riordan graphs}\label{Riordan-matrix-sec}
 Let $\kappa[[z]]$ be the ring of formal power series over an
integral domain $\kappa$.
 If there exists a pair of generating functions $(g,f)\in \kappa[[z]]\times \kappa[[z]]$, $f(0)=0$ such
 that for $j\ge 0$,
\begin{eqnarray*}
 g f^j=\sum_{i\ge0}\ell_{i,j}z^i,
\end{eqnarray*}
then the matrix  $L=[\ell_{i,j}]_{i,j\ge0}$ is called a  {\it
Riordan matrix} (or, a {\it Riordan array})  over $\kappa$ generated
by $g$ and $f$. Usually, we write $L=(g,f)$. Since $f(0)=0$, every
Riordan matrix $(g,f)$ is an infinite lower triangular matrix. If
a Riordan matrix is invertible, it is called {\it proper}. Note that
$(g,f)$ is invertible if and only if $g(0)\ne0$, $f(0)=0$ and
$f^\prime(0)\ne0$.

For a Riordan matrix $(g,f)$ over $\mathbb{Z}$, the matrix
$L=[\ell_{i,j}]_{i,j\ge0}$ defined by
\begin{eqnarray*}
 \ell_{i,j}\equiv[z^i]gf^j\;({\text{mod}}\;2),
\end{eqnarray*}
is called a {\it binary Riordan matrix}, and it is denoted by ${\cal
B}(g,f)$. The {\em leading principal matrix of order $n$} in ${\cal
B}(g,f)$ is denoted by ${\cal B}(g, f)_n$.


The following result is well known as the {\em Fundamental
Theorem of Riordan matrices} ({\em FTRM}).
 \begin{lemma}[FTRM, \cite{shap}] \label{FTRM-lem}
 Let $R=(g,f)$ be a Riordan matrix and let $R\,A=B$ where $A$ and $B$ are two infinite vectors with the generating functions $a(z)$ and $b(z)$,
 respectively. Then $b=g\, a(f)$.
 \end{lemma}



The following definition gives the notion of a Riordan graph in both {\em labelled} and {\em unlabelled} cases. We note that throughout this paper the graphs are assumed to be {\em labelled} unless otherwise specified.

 \begin{defin}\label{main-labelled} {\rm A simple {\em
labelled} graph $G$ on the vertex set $[n]$ is a {\em Riordan graph}
of order $n$ if the adjacency matrix of $G$, ${\cal A}(G)$, is an
$n\times n$ symmetric (0,1)-matrix such that
 \begin{itemize}
\item its main diagonal entries are all 0, and
\item its lower triangular part below the main diagonal is ${\cal B}(g,f)_{n-1}$
\end{itemize}
 for some Riordan matrix $(g,f)$ over $\mathbb{Z}$, i.e., by using Riordan language,
\begin{eqnarray}
{\cal A}(G)={\cal B}(zg,f)_n+{\cal B}(zg, f)_n^T.\label{def1}
 \end{eqnarray}
 We denote such $G$
  by $G_n(g,f)$, or simply by $G_n$ when the pair $(g,f)$ is understood from the context, or it is not important.
A simple {\em unlabelled} graph is a {\em Riordan graph} if at least
one of its labelled copies is a Riordan graph.}
\end{defin}

We note that the choice of the functions $g$ and $f$ in Definition~\ref{main-labelled} may not be unique. If $G$ is a Riordan graph and ${\cal A}(G)=[r_{ij}]_{i,j\ge1}$, then for $i>j\geq 1$,
\begin{eqnarray}\label{def2}
r_{i,j}\equiv[z^{i-1}]zgf^{j-1}\;({\text{mod}}\;2)\equiv[z^{i-2}]gf^{j-1}\;({\text{mod}}\;2).
 \end{eqnarray}

For example, the Riordan graph  $G_n\left({1\over 1-z},{z\over
1-z}\right)$ is called the {\it Pascal graph} of order $n$, and is
denoted by $PG_n$. For another example, the {\em Catalan graph}
$CG_6(C,zC)$ shown below, where
$$C={1-\sqrt{1-4z}\over 2z}=\sum_{n\geq 0}\frac{1}{n+1}{2n\choose
n}z^n=1+z+2z^2+5z^3+14z^4+\cdots,$$ is given by the adjacency matrix
$\mathcal{A}(CG_6)$:
  \begin{center}
  \begin{tabular}{cc}
 $\mathcal{A}(CG_6)=\left[\begin{array}{cccccc}
        0 & 1 & 1 & 0 & 1 & 0 \\
        1 & 0 & 1 & 0 & 1 & 0 \\
        1 & 1 & 0 & 1 & 1 & 1 \\
        0 & 0 & 1 & 0 & 1 & 0 \\
        1 & 1 & 1 & 1 & 0 & 1 \\
        0 & 0 & 1 & 0 & 1 & 0
      \end{array}\right]$& \scalebox{1} 
 {
 \begin{pspicture}(0,0)(3.02625,0.0939063)
 \psdots[dotsize=0.12](0.08,1.2104688)
 \psdots[dotsize=0.12](1.48,1.2104688)
 \psdots[dotsize=0.12](2.88,1.2104688)
 \psdots[dotsize=0.12](0.38,-0.78953123)
 \psdots[dotsize=0.12](1.48,-1.3895313)
 \psdots[dotsize=0.12](2.58,-0.78953123)
 \usefont{T1}{ptm}{m}{n}
 \rput(0.07859375,1.5204687){2}
 \usefont{T1}{ptm}{m}{n}
 \rput(1.5009375,1.5204687){4}
 \usefont{T1}{ptm}{m}{n}
 \rput(2.8953125,1.5204687){6}
 \usefont{T1}{ptm}{m}{n}
 \rput(0.166875,-0.93953127){1}
 \usefont{T1}{ptm}{m}{n}
 \rput(1.3076563,-1.5395312){3}
 \usefont{T1}{ptm}{m}{n}
 \rput(2.6695313,-0.95953125){5}
 \psline[linewidth=0.04cm](0.08,1.2104688)(0.38,-0.78953123)
 \psline[linewidth=0.04cm](0.08,1.2304688)(1.48,-1.3695313)
 \psline[linewidth=0.04cm](0.08,1.2104688)(2.58,-0.78953123)
 \psline[linewidth=0.04cm](1.48,1.2104688)(1.48,-1.3495313)
 \psline[linewidth=0.04cm](1.48,1.2304688)(2.58,-0.74953127)
 \psline[linewidth=0.04cm](2.86,1.2104688)(1.48,-1.4095312)
 \psline[linewidth=0.04cm](2.88,1.2104688)(2.6,-0.76953125)
 \psline[linewidth=0.04cm](0.38,-0.78953123)(2.58,-0.78953123)
 \psline[linewidth=0.04cm](0.38,-0.78953123)(1.48,-1.3895313)
 \psline[linewidth=0.04cm](1.48,-1.3895313)(2.58,-0.78953123)
 \end{pspicture}
 }
  \end{tabular}
  \end{center}

\ \\[2mm]

\begin{defin}\label{def-prop-graph} {\rm A Riordan graph $G_n(g,f)$ is {\em proper} if the binary Riordan matrix ${\cal B}(g,f)_{n-1}$ is proper.
Thus, in a proper Riordan graph a vertex $i$ is adjacent to the
vertex $i+1$ for $i\ge1$.} \end{defin}

If a Riordan graph $G_n(g,f)$ is proper then
the Riordan matrix $(g,f)$ is also proper because $g(0)\;({\rm mod}\ 2)\equiv f'(0)\;({\rm mod}\ 2)\equiv 1$. The converse to this statement is not true. For
instance, $(1,2z+z^2)$ is a proper Riordan matrix but
$G_n(1,2z+z^2)$ is not a proper Riordan graph.

The following theorem about the adjacency matrices of Riordan graphs is a key result in \cite{CJKM}.

\begin{theorem}[Riordan Graph Decomposition, \cite{CJKM}]\label{e:th}  Let $G_n=G_n(g,f)$ be a Riordan graph  with $[z^1]f=1$. Then
\begin{itemize}
\item[{\rm(i)}] The adjacency matrix $\mathcal{A}(G_{n})$ satisfies
\begin{eqnarray}\label{e:bm}
\mathcal{A}(G_{n})=P^{T}\left[
\begin{array}{cc}
X & B \\
B^{T} & Y
\end{array}
\right]P
\end{eqnarray}
where $P=\left[e_{1}\;|\; e_{3}\;|\; \cdots\;|\; e_{2\lceil
n/2\rceil -1}\;|\; e_{2}\;|\; e_{4} \;|\; \cdots \;|\; e_{2\lfloor
n/2\rfloor }\right] ^{T}$ is the $n\times n$ permutation matrix and $e_{i}$ is
the elementary column vector with the $i$th entry being $1$ and
the others entries being~$0$.
\item[{\rm(ii)}] The matrix $X$ is the adjacency matrix of the
induced subgraph of $G_n(g,f)$ by the odd indexed vertex set
$V_o=\{2i-1\;|\;1\le i\le \lceil n/2\rceil\}$. In particular, the
induced subgraph $\left<V_o\right>$ is a Riordan graph of order
${\lceil n/2\rceil}$ given by $G_{\lceil
n/2\rceil}(g^{\prime}(\sqrt{z}) ,f(z))$.
\item[{\rm(iii)}] The matrix $Y$ is the adjacency matrix of the induced
subgraph of $G_n(g,f)$ by the even indexed vertex set
$V_e=\{2i\;|\;1\le i\le \lfloor n/2\rfloor\}$. In particular, the
induced subgraph $\left<V_e\right>$ is a Riordan graph of order
${\lfloor n/2\rfloor }$ given by $G_{\lfloor n/2\rfloor
}\left(\left(\frac{gf}{z}\right)^{\prime}(\sqrt{z}),f(z)\right)$.
\item[{\rm(iv)}] The matrix $B$ representing the edges between $V_o$ and $V_e$
can be expressed as the sum of binary Riordan matrices as follows:
\begin{align*}
B={\cal B}(z\cdot(gf)^{\prime }(\sqrt{z}),f(z))_{\lceil n/2\rceil
\times\lfloor n/2\rfloor}+{\cal B}((zg)^{\prime}(\sqrt{z}),f(z))
_{\lfloor n/2\rfloor\times\lceil n/2\rceil }^{T}.
\end{align*}
\end{itemize}
\end{theorem}

\subsection{Families of Riordan graphs}

There are many families of Riordan graphs introduced in~\cite{CJKM}. Below we list those of them most relevant to this paper, along with examples of subfamilies. \\[-2mm]

\noindent
{\bf Riordan graphs of the {\em Appell type}.} This class of graphs is defined by an {\em Appell matrix} $(g,z)$, and thus it is precisely the class of  Toeplitz graphs. Examples of graphs in this class are
\begin{itemize}
\item the null graphs $N_n$ defined by $G_n(0,z)$;
\item the path graphs $P_n$ defined by $G_n(1,z)$;
\item the complete graphs $K_n$ defined by $G_n\left(\frac{1}{1-z},z\right)$; and
\item the complete bipartite graphs $K_{\lfloor{n\over2}\rfloor, \lceil{n\over2}\rceil}$ defined by $G_n\left({1\over 1-z^2},z\right)$.
\end{itemize}

\noindent
{\bf Riordan graphs of the {\em Bell type}.} This class of graphs is defined by a {\em Bell matrix} $(g,zg)$. Examples of graphs in this class are
\begin{itemize}
\item the null graphs $N_n$ defined by $G_n(0,0)$;
\item the path graphs $P_n$ defined by $G_n(1,z)$;
\item the Pascal graphs $PG_n$ defined by $G_n\left({1\over 1-z},{z\over
1-z}\right)$; and
\item the {\em Catalan graphs} $CG_n$ defined by $G_n\left({1-\sqrt{1-4z}\over 2z},{1-\sqrt{1-4z}\over
2}\right)$.
\end{itemize}

\noindent
{\bf Riordan graphs of the {\em checkerboard type.}} This class of graphs is defined by a {\em checkerboard matrix} $(g,f)$ such that $g$ is an even function and $f$ is an odd function. Examples of graphs in this class are
\begin{itemize}
\item the null graphs $N_n$ defined by $G_n(0,c)$ for any constant $c$;
\item the path graphs $P_n$ defined by $G_n(1,z)$; and
\item the complete bipartite graphs $K_{\lfloor{n\over2}\rfloor, \lceil{n\over2}\rceil}$ defined by $G_n\left({1\over 1-z^2},z\right)$.
\end{itemize}

\noindent
{\bf Riordan graphs of the {\em derivative type.}} This class of graphs is defined by functions $(f',f)$. Examples of graphs in this class are
\begin{itemize}
\item the null graphs $N_n$ defined by $G_n(0,c)$ for any constant $c$; and
\item the path graphs $P_n$ defined by $G_n(1,z)$.
\end{itemize}

\noindent
{\bf o-decomposable Riordan graphs} standing for {\em odd decomposable Riordan graphs}. This class of graphs is defined by requiring in (\ref{e:bm}) $Y=O$, where $O$ is the zero matrix of $Y$'s size.\\[-2mm]

\noindent
{\bf e-decomposable Riordan graphs} standing for {\em even decomposable Riordan graphs}. This class of graphs is defined by requiring in (\ref{e:bm}) $X=O$. \\[-2mm]

\subsection{Several known results on Riordan graphs}

 The following lemma is a corollary of Theorem~\ref{e:th}.
\begin{lemma}[\cite{CJKM}]\label{lem12}
Every adjacency matrix ${\cal A}={\cal A}_{n}(g,zg)$ of a Riordan
graph $G_{n}(g,zg)$ of the Bell type can be expressed as the block
matrix
\begin{eqnarray}\label{e:bm1}
P{\cal A}P^{T}=\left[
\begin{array}{cc}
X & B \\
B^{T} & O%
\end{array}
\right],
\end{eqnarray}
where $$X={\cal A}_{\lceil n/2\rceil }(g^{\prime
}(\sqrt{z}),zg(z))$$ and
$$B\equiv {\cal B}(zg(z),zg(z))_{\lceil n/2\rceil, \lfloor n/2\rfloor}+{\cal B}((zg)^{\prime }(\sqrt{z}),zg(z)) _{\lfloor n/2\rfloor,\lceil n/2\rceil}^T.$$ Thus, any Riordan graph of the Bell type is o-decomposable.
 \end{lemma}

The following lemma is given by Theorem~3.6 (iv) in~\cite{CJKM}.

 \begin{lemma}[\cite{CJKM}]\label{check} A graph $G_n$ is a Riordan graph of the checkerboard type if and only if $G_n$ is bipartite with bipartitions $V_o$ and $V_e$.
 \end{lemma}

In this paper, we also need the following lemmas.

\begin{lemma}[\cite{SG}]\label{PC0}  For the Pascal graph $PG_n$,  we have
 $m(PG_{2^{k}})=3^{k}-2^{k}$ and $m( PG_{2^{k}+1})=3^{k}$.  \end{lemma}

\begin{lemma}[\cite{CJKM}]\label{PC}  For the Catalan graph $CG_n$, we have $m\left(CG_{2^{k}}\right)=\frac{3^{k}-1}{2}$ and $m\left(CG_{2^{k}+1}\right)=\frac{3^{k}-1}{2}+2^{k}$.
 \end{lemma}

\begin{defin}[\cite{CJKM}]\label{ioie-decomposable-def} {\rm Let
$G_n=G_{n}(g,f)$ be a proper Riordan graph
 with the odd and even vertex sets $V_{o}$ and $V_{e}$, respectively.
\begin{itemize}
\item If $\left< V_{o}\right> \cong G_{\lceil n/2\rceil}(g,f)$ and $\left<
V_{e}\right>$ is a null graph then $G_n$ is {\em io-decomposable}.
\item If $\left< V_{o}\right>$ is a null graph and $\left< V_{e}\right> \cong
G_{\lfloor n/2\rfloor}(g,f)$ then $G_n$ is {\em ie-decomposable}.
\end{itemize}
``io'' and ``ie'' stand for ``isomorphically odd'' and
``isomorphically even'', respectively.}\end{defin}

 \begin{lemma}[\cite{CJKM}]\label{lm2}
Let $G_n=G_n(g,zg)$ be an io-decomposable Riordan graph. Then, we have
 $$\chi(G_n)=\omega(G_n)=\lceil\log_2 n\rceil+1,$$
 where $\chi(G_n)$ and $\omega(G_n)$ are the chromatic number and clique number of $G_n$, respectively.
\end{lemma}

\begin{lemma}[\cite{CJKM}]\label{e:lem3}
A Riordan graph $G_{n}(g,zg)$ is io-decomposable if and only if
\begin{eqnarray*}
g^{2}\equiv g^{\prime}, \text{ i.e.\ $[z^j]g\equiv[z^{2j+1}]g$.}
\end{eqnarray*}
\end{lemma}

\subsection{The spectral graph theory}

Let $G$ be a graph and $\mathcal{A}(G)$ be its adjacency matrix.
The eigenvalues $\lambda_1\geq \cdots \geq \lambda_n$ of
$\mathcal{A}(G)$ are said to be the {\em eigenvalues} of $G$, and
the eigenvalues form the {\em spectrum} of this graph. The {\em
largest} eigenvalue $\lambda_1$ is called the {\it spectral radius}
of $G$.  The determinant of $G$, denoted by $\det(G)$, is the
determinant of $\mathcal{A}(G)$. The number of positive and negative
eigenvalues of a graph are called {\em positive} and {\em negative
inertias} of the graph, respectively. We denote the positive and
negative inertias of a Riordan graph $G_n$ by $n^+(G_n)$ and
$n^-(G_n)$, respectively. If $G_n$ is understood from the context,
we simply write $n^+$ and $n^-$ for the inertias.

  Let $M$ be an $i\times j$ matrix. The {\em null space} of $M$ is the set of all $j$-dimensional column vectors $X$ such that $M X=0$. The {\em dimension} of the null space of $M$ is called the {\em nullity} of the matrix $M$, and
is denoted by $\eta(M)$, or just by $\eta$ when $M$ is understood
from the context. The nullity of $G$ is
$\eta(G)=\eta(\mathcal{A}(G))$. A graph $G$ is called {\em
singular} if its adjacency matrix is singular.

\begin{defin} {\rm The {\em Laplacian matrix} $L(G)$ and {\em signless Laplacian matrix} $Q(G)$ of a graph $G$ are
 defined, respectively, as $D(G)-\mathcal{A}(G)$ and $D(G)+\mathcal{A}(G)$,  where $D(G)$ is the diagonal matrix of vertex
 degrees of $G$.} \end{defin}

\begin{defin} {\rm The {\em Laplacian spectrum} of a graph $G$ is the sequence of its Laplacian eigenvalues,
i.e.\  $\mu_1\ge \mu_2\ge \cdots\ge\mu_{n-1}\ge\mu_n=0$. We let $Lspec(G)$ denote the set $\{\mu_1,\mu_2,\ldots,\mu_n\}$.
The {\em signless Laplacian spectrum} of $G$ is the sequence of its signless Laplacian eigenvalues,
i.e. $q_1\ge q_2 \ge \cdots\ge q_n$.} \end{defin}

It is well known that the average degree $\frac{2m}{n}$ of a graph $G$ with $n$ vertices and $m$ edges is a lower bound for the spectral radius $\lambda_1(G)$, i.e.
 \begin{equation}
 \lambda_1(G)\ge \frac{2m}{n}.  \label{av}
 \end{equation}
 Moreover, for the signless Laplacian spectral radius $q_1(G)$, we have
 \begin{equation}
 q_1(G)\ge 2\lambda(G)\ge \frac{4m}{n}.  \label{av2}
 \end{equation}

 The following lemma gives a relation between the eigenvalues of a real symmetric matrix and the eigenvalues of its {\em partitioned matrix}.

 \begin{lemma}[\cite{cv}]\label{lem31}
 Let $A$ be a real symmetric matrix with eigenvalues $\lambda_1\ge \lambda_2\ge \cdots \ge \lambda_n$. Given a partition $\{1,2,\ldots,n\}=V_1\mathbin{\dot{\cup}}V_2 \mathbin{\dot{\cup}} \cdots \mathbin{\dot{\cup}} V_k$
with $|V_i|=n_i>0$, consider the corresponding blocking
$A=(A_{i,j})$, where $A_{i,j}$ is an $n_i\times n_j$ block and $1\le
i,j \le k$. Let $e_{i,j}$ be the sum of the entries in $A_{i,j}$ and
set the matrix $B:=(e_{i,j}/n_i)$ for $1\le i,j \le k$. Then the
eigenvalues of $B$ interlace those of $A$, i.e. $\lambda_i\ge \rho_i
\ge \lambda_{n-k-i}$ for $1\le i \le k$, where $\rho_i$ is the $i$th
largest eigenvalue of $B$.
 \end{lemma}

The following lemma gives a relation between the adjacency matrix of
a graph $G$ and the {\em clique number} $\omega(G)$ of $G$, which is the size of a maximal clique in $G$.

\begin{lemma}[\cite{cv}]\label{lm1}
If $G$ is a graph with the adjacency matrix $A$ and the clique
number $\omega(G)$, then $$\max{\{X^T A X\,;\, X\ge0,\, {\bf 1}^T
X=1\}}=1-\frac{1}{\omega(G)},$$ where ${\bf 1}$ is the vector of the proper size that contains only $1$s .
\end{lemma}

In this paper, we also need the following lemmas.

\begin{lemma}[\cite{cv}] \label{par}
Let $M$ be a symmetric $n\times n$ matrix with real entries. If
$M=\left[
                                                                     \begin{array}{cc}
                                                                       P & Q \\
                                                                       Q^T & R \\
                                                                     \end{array}
                                                                   \right]$, then $\lambda_1(M)+\lambda_n(M)\le \lambda_1(P)+\lambda_1(R)$, where $\lambda_1(M)$ and $\lambda_n(M)$ stand for the largest and smallest eigenvalies of $M$, respectively.
\end{lemma}

\begin{lemma}[\cite{Ber}]\label{lem41}
Let $A$ and $B$ be real matrices of orders $n\times n$ and $n\times
m$, respectively, and let $r={\rm rank}(B)$ and ${\cal
A}:=\left[
                                                                                                                                                    \begin{array}{cc}
                                                                                                                                                      A & B \\
                                                                                                                                                      B^T & O \\
                                                                                                                                                    \end{array}
                                                                                                                                                  \right]
$. Then $n^-({\cal A})\ge r$, $\eta({\cal A})\ge 0$ and $n^+({\cal
A})\ge r$. If, in addition, $n=m$ and $B$ is nonsingular, then
$n^-({\cal A})=n^+({\cal A})=n$ and $\eta({\cal A})=0$.
\end{lemma}

\begin{lemma}[\cite{Ber}]\label{lem61}
Let $A\in {\mathbb{R}}^{n \times n}$, $B\in {\mathbb{R}}^{m \times
m}$, $C\in {\mathbb{R}}^{n \times m}$, and $A$ and $B$
be symmetric. Also, let ${\cal A}_0:=\left[
                                                                                                                                                                               \begin{array}{cc}
                                                                                                                                                                                 A & O \\
                                                                                                                                                                                 O & B\\
                                                                                                                                                                               \end{array}
                                                                                                                                                                             \right]
$, ${\cal A}:=\left[
                                                                                                                                                                               \begin{array}{cc}
                                                                                                                                                                                 A & C\\
                                                                                                                                                                                 C^T & B\\
                                                                                                                                                                               \end{array}
                                                                                                                                                                             \right]
$, $\sigma_{\max}(C)$ be the largest singular value of $C$,
and $\eta:=\min\limits_{\substack{i=1,\ldots,n{}\\j=1,\ldots,m}}
|\lambda_i(A)-\lambda_j(B)|$. Then, for all $i=1,\ldots,n+m,$
$$|\lambda_i({\cal A})-\lambda_i({\cal A}_0)|\le \frac{2\sigma_{\max}^2 (C)}{\eta+\sqrt{\eta^2+4\sigma_{\max}(C)}}.$$
\end{lemma}

\begin{lemma}[\cite{cv,HT}]\label{lm3111}
For any graph $G$ with chromatic number $\chi(G)$ and clique number $\omega(G)$, we have
 \begin{itemize}
\item[{\rm(i)}]  $\chi(G)-1\le \lambda_1(G)\le n\,\left(1-\frac{1}{\omega(G)}\right)$.
\item[{\rm(ii)}]  $\chi(G)\ge 1+\frac{\lambda_1(G)}{|\lambda_n(G)|}$.
\item[{\rm(iii)}]  $\chi(G)\ge 1+\frac{\lambda_1(G)}{\mu_1(G)-\lambda_1(G)}$.
\end{itemize}
\end{lemma}

 \begin{lemma}[\cite{MERIS}]\label{lem51}
 Let $G$ be a graph of order $n$. Then  $\mu_1(G)\leq n$, where the equality holds if and only if the complement $\overline{G}$ of $G$ is
 disconnected. Moreover, $Lspec(\overline{G})=\{n-\mu_{n-1},n-\mu_{n-2},\dots,n-\mu_{1},0\}$.
 \end{lemma}

 \begin{lemma}[\cite{FI}]\label{mb1} Let $G$ $(\ncong K_n)$ be a graph of order $n$. Then $\mu_{n-1}(G)\leq \delta(G)$, where $\delta(G)$ is the minimum vertex degree in $G$.
 \end{lemma}

\begin{lemma}[\cite{MERIS}]\label{max1}
Let $G$ be a graph on $n$ vertices with at least one edge. Then
$\mu_1(G)\ge \Delta(G)+1$, where $\Delta(G)$ is the maximum vertex degree in $G$. Moreover, if $G$ is connected, then the
equality holds if and only if $\Delta(G)=n-1$.
\end{lemma}

The following result states for which $n$ the Catalan number
$C_n$ is odd.
\begin{lemma}[\cite{DS}]\label{lem4}
 The only Catalan numbers $C_n$ that are odd are those for which $n=2^k-1$  for $k\ge 0$.
 In particular,  $C_{2i}\equiv 0$ for $i\in \mathbb{N}$.
 \end{lemma}

 \section{Eigenvalues of Riordan graphs}\label{eigenvalues-sec}

 In this section we present several lower and upper bounds on the eigenvalues and the (signless) Laplacian eigenvalues of Riordan graphs, which will improve  \eqref{av} and \eqref{av2} for this class of graphs.
 To achieve our goals, we begin with a general result valid for all graphs.

 \begin{theorem}\label{thm9}
 Let $G=(V, E)$ be a graph with $|V|=n$ and $|E|=m$, and let $W$ be a nonempty subset of $V$ with $|W|=k$. Then
 \begin{eqnarray}
 &{\rm(i)}&   \lambda_1(G)   \ge  \frac{m_1}{k}+ \frac{m_2}{n-k} + \sqrt{\left( \frac{m_1}{k}-\frac{m_2}{n-k}\right)^2+\frac{(m-m_1-m_2)^2}{k(n-k)}}.  \nonumber\\[2mm]
 &{\rm(ii)}&   \mu_1(G)     \ge   \frac{n (m-m_1-m_2)}{k (n-k)}. \nonumber\\[2mm]
 &{\rm(iii)}&    q_1(G)       \ge  \frac{1}{2} \left( c_1+c_2+\sqrt{\Big(c_1-c_2\Big)^2+\frac{4(m-m_1-m_2)^2}{k(n-k)} }  \right).\nonumber\\[2mm]
 &{\rm(iv)}&   q_2(G)         \ge  \frac{1}{2} \left( c_1+c_2-\sqrt{\Big(c_1-c_2\Big)^2+\frac{4(m-m_1-m_2)^2}{k(n-k)} }  \right) \ge q_n(G)\nonumber
 \end{eqnarray}   where $c_1=\frac{m+3m_1-m_2}{k}$, $c_2=\frac{m-m_1+3m_2}{n-k}$, $m_1=m(\langle W \rangle)$ and $m_2=m(\langle V\backslash W \rangle)$.
 \end{theorem}

 \begin{pf} Let $ \mathcal{W}$ and $\mathcal{U}$ be, respectively, the adjacency matrices of the induced subgraphs $\langle W \rangle$ and  $\langle V\backslash W \rangle$. Let $\mathcal{A}$ be the adjacency matrix of $G$. Then there exists a permutation matrix $P$ such that
 $$P\mathcal{A}P^T=\left[
           \begin{array}{cc}
             \mathcal{W} & \mathcal{F} \\
             \mathcal{F}^T &  \mathcal{U} \\
           \end{array}
         \right].
 $$ Using Lemma~\ref{lem31} on $P\mathcal{A}P^T$, we obtain
 $$B_1:=\left[
           \begin{array}{cc}
            \frac{2m_1}{k_1} & \frac{m_3}{k_1} \\
                             &                  \\
             \frac{m_3}{k_2} &  \frac{2m_2}{k_2} \\
           \end{array}
         \right]
 $$ where $k_1:=|W|$ and $k_2:=n-k_1=|V\backslash W|$, and $m_3$ stands for the number of 1s in the submatrix $\mathcal{F}$, which is the number of edges between $W$ and $V\backslash W$.
 It is not difficult to see that the eigenvalues of $B_1$ are given by
 $$\lambda_1(B_1)=\frac{m_1}{k_1}+ \frac{m_2}{k_2} + \sqrt{\left( \frac{m_1}{k_1}-\frac{m_2}{k_2}\right)^2+\frac{m_3^2}{k_1\,k_2}}$$ and
 $$\lambda_2(B_1)=\frac{m_1}{k_1}+ \frac{m_2}{k_2} - \sqrt{\left( \frac{m_1}{k_1}-\frac{m_2}{k_2}\right)^2+\frac{m_3^2}{k_1\,k_2}}.$$ By Lemma \ref{lem31}, we have $\lambda_1(G)\ge \lambda_1(B_1) $. This completes the proof of (i).

  Suppose that $L=D-\mathcal{A}$ is the Laplacian matrix of $G$. Then by the same permutation matrix $P$ as above, we obtain
   $$PLP^T=PDP^T-P\mathcal{A}P^T=\left[
           \begin{array}{cc}
             D_\mathcal{W} & O \\
             O^T &  D_\mathcal{U} \\
           \end{array}
         \right]-
   \left[
           \begin{array}{cc}
             \mathcal{W} & \mathcal{F} \\
             \mathcal{F}^T &  \mathcal{U} \\
           \end{array}
         \right]=
            \left[
           \begin{array}{cc}
              D_\mathcal{W}-\mathcal{W} & -\mathcal{F} \\
             -\mathcal{F}^T &  D_\mathcal{U}-\mathcal{U} \\
           \end{array}
         \right]
 $$ where $ D_\mathcal{W}$ and $D_\mathcal{U}$ are, respectively, the degree matrices with the same row indices as $\mathcal{W}$ and $\mathcal{U}$, and $O$ is the zero matrix of the proper size.
 Again, using Lemma \ref{lem31} on $PLP^T$, we obtain
 $$B_2:=\left[
           \begin{array}{cc}
            \frac{m_3}{k_1} & -\frac{m_3}{k_1} \\
                             &                  \\
             -\frac{m_3}{k_2} &  \frac{m_3}{k_2} \\
           \end{array}
         \right].$$
 It is easy to see that the eigenvalues of $B_2$ are given by
 $$\lambda_1(B_2)= \frac{nm_3}{k_1 k_2}  ~~~\mbox{and}~~~ \lambda_2(B_2)=0.$$ By Lemma \ref{lem31}, we have $\mu_1(G)\ge \lambda_1(B_2)$. Hence, (ii) is proved.

    Suppose that $Q=D+A$ is the signless Laplacian matrix of $G$. Then, by the same permutation matrix $P$ as above, we obtain
   $$PQP^T=PDP^T+P\mathcal{A}P^T=\left[
           \begin{array}{cc}
             D_\mathcal{W} & O \\
             O^T &  D_\mathcal{U} \\
           \end{array}
         \right]+
   \left[
           \begin{array}{cc}
             \mathcal{W} & \mathcal{F} \\
             \mathcal{F}^T &  \mathcal{U} \\
           \end{array}
         \right]=
            \left[
           \begin{array}{cc}
              D_\mathcal{W}+\mathcal{W} & \mathcal{F} \\
             \mathcal{F}^T &  D_\mathcal{U}+\mathcal{U} \\
           \end{array}
         \right]. $$
         Using Lemma~\ref{lem31} on $PQP^T$, we obtain
 $$B_3=\left[
           \begin{array}{cc}
            \frac{4m_1+m_3}{k_1} & \frac{m_3}{k_1} \\
                             &                  \\
             \frac{m_3}{k_2} &  \frac{4m_2+m_3}{k_2} \\
           \end{array}
         \right].$$  The eigenvalues of $B_3$ are given by
 $$\lambda_1(B_3)= \frac{1}{2} \left( \frac{4m_1+m_3}{k_1}+ \frac{4m_2+m_3}{k_2}+\sqrt{\Big( \frac{4m_1+m_3}{k_1}-\frac{4m_2+m_3}{k_2}\Big)^2+\frac{4m_3^2}{k_1 k_2} }  \right)$$ and
 $$\lambda_2(B_3)= \frac{1}{2} \left( \frac{4m_1+m_3}{k_1}+\frac{4m_2+m_3}{k_2}-\sqrt{\Big( \frac{4m_1+m_3}{k_1}-\frac{4m_2+m_3}{k_2}\Big)^2+\frac{4m_3^2}{k_1 k_2} }  \right).$$  One can see that $\lambda_2(B_3)>0$.
  By Lemma \ref{lem31}, we have $q_1(G)\ge \lambda_1(B_3) $ and $q_2(G) \ge \lambda_2(B_3) \ge q_n(G)$.
 This completes the proofs for  (iii) and (iv).
 \end{pf}

 \begin{corollary}\label{cor-thm9}
 Let $G=(V, E)$ be a graph with $|V|=n$ and $|E|=m$, and let $W$ be a nonempty subset of $V$ with $|W|=k$. Also, let $m_1=m(\langle W \rangle)$ and $m_2=m(\langle V\backslash W \rangle)$, and $m> m_1+m_2+2\sqrt{m_1\,m_2}$. Then $\lambda_n(G)<0$.
\end{corollary}

\begin{pf}
For any graph $G$, $\lambda_2(B_1) \ge \lambda_n(G)$, where $B_1$ and $\lambda_2(B_1)$ can be found in the proof of Theorem~\ref{thm9}. Since $m> m_1+m_2+2\sqrt{m_1\,m_2}$, we have
$$ \lambda_n(G)  \le   \frac{m_1}{k}+ \frac{m_2}{n-k} - \sqrt{\left( \frac{m_1}{k}-\frac{m_2}{n-k}\right)^2+\frac{(m-m_1-m_2)^2}{k(n-k)}}< 0.$$
\end{pf}

 The following theorem, the main result in this section, follows from  Theorems~\ref{e:th} and~\ref{thm9}.
 \begin{theorem} \label{co1}
 Let $G_n=G_n(g, f)$ be a Riordan graph with $n$ vertices and $m$ edges. Then
 \begin{eqnarray}
 &{\rm(i)}&  \lambda_1(G_n)  \ge  \frac{m_1}{n_1}+ \frac{m_2}{n-n_1} + \sqrt{\left( \frac{m_1}{n_1}-\frac{m_2}{n-n_1}\right)^2+\frac{(m-m_1-m_2)^2}{n_1(n-n_1)}}.  \label{lam}\\[2mm]
 &{\rm(ii)}&\mu_1(G_n)     \ge   \frac{n (m-m_1-m_2)}{n_1 (n-n_1)}. \nonumber\\[2mm]
 &{\rm(iii)}& q_1(G_n)      \ge  \frac{1}{2} \left( c_1+c_2+\sqrt{\Big(c_1-c_2\Big)^2+\frac{4(m-m_1-m_2)^2}{n_1(n-n_1)} }  \right).\label{lam2}\\[2mm]
 &{\rm(iv)}&  q_2(G_n)      \ge  \frac{1}{2} \left( c_1+c_2-\sqrt{\Big(c_1-c_2\Big)^2+\frac{4(m-m_1-m_2)^2}{n_1(n-n_1)} }  \right) \ge q_n(G),\nonumber
  \end{eqnarray} where $n_1=\lceil\frac{n}{2} \rceil$, $m_1=m(G_{n_1}(g'(\sqrt{z}), f))$, $m_2=m(G_{n-n_1}((\frac{gf}{z})'(\sqrt{z}), f))$,  $c_1=\frac{m+3m_1-m_2}{n_1}$ and $c_2=\frac{m-m_1+3m_2}{n-n_1}$.
  \end{theorem}

 \begin{proposition}\label{prop-improvement}
 For a Riordan graph $G_n$, the lower bound in (\ref{lam}) is greater than or equal to the lower bound in~(\ref{av}). \end{proposition}
 \begin{pf}
 Suppose that $n=2k$, and $m_3$ is defined as in the proof of Theorem~\ref{thm9}. We have
 $$(m_1-m_2)^2\ge 0~~\Rightarrow ~(m_1+m_2)^2\ge 4m_1 m_2~\Rightarrow$$
 $$m_3^2\ge m_3^2+(m_1+m_2)^2+2m_3(m_1+m_2)+4m_1m_2-2(m_1+m_2)^2-2m_3(m_1+m_2)$$
 $$\Rightarrow ~m_3^2\ge m^2+4m_1m_2-2m m_1-2m m_2~\Rightarrow ~$$
 $$\frac{(m_1-m_2)^2+m_3^2}{k^2}\ge \frac{m^2+(m_1+m_2)^2-2m(m_1+m_2)}{k^2} ~\Rightarrow ~$$
 $$\frac{m_1+m_2}{k}+\sqrt{\frac{(m_1-m_2)^2+m_3^2}{k^2}}\ge \frac{m}{k}~\Rightarrow $$
 $$~\frac{m_1}{n_1}+ \frac{m_2}{n_2} + \sqrt{\left( \frac{m_1}{n_1}-\frac{m_2}{n_2}\right)^2+\frac{m_3^2}{n_1n_2}}\ge \frac{2m}{n}$$
 as $n_1=n_2=k$. Similarly, the result for $n=2k+1$ follows.  \end{pf}

\begin{proposition}\label{prop-improvement1}
 For a Riordan graph $G_n$, the lower bound in (\ref{lam2}) is greater than or equal to the lower bound in~(\ref{av2}). \end{proposition}
 \begin{pf}
Since $q_1\geq 2\lambda_1$ for any graph, the result follows from Proposition~\ref{prop-improvement}.
\end{pf}

The following results are obtained by Theorems~\ref{e:th} and~\ref{co1}.

 \begin{corollary}\label{cor-App-t-RG}
 Let $G_n=G_n(g, z)$ be a Riordan graph of the Appell type with $n$ vertices and $m$ edges. Let $n$ be even and $m_1=m(G_{n/2}(g'(\sqrt{z}), z))$. Then
  \begin{itemize}
\item[{\rm(i)}] The lower bounds in (\ref{lam}) and (\ref{av}) coincide for the spectral radius $\lambda_1(G_n)$.
\item[{\rm(ii)}] The lower bounds in (\ref{lam2}) and (\ref{av2}) coincide for the signless Laplacian spectral radius $q_1(G_n)$.
\item[{\rm(iii)}]  $\mu_1(G_n)   \ge   \frac{4 (m-2m_1)}{n}$.
\item[{\rm(iv)}]   $q_2(G_n)     \ge  \frac{8m_1}{n}   \ge q_n(G)$.
\item[{\rm(v)}]  If $m> 4m_1$, then   $ \lambda_n(G_n) \le \frac{8m_1-2m}{n}<0$.
 \end{itemize}
  \end{corollary}

 \begin{corollary}\label{cor2}
 Let $G_n=G_n(g, f)$ be a Riordan graph with $n$ vertices and $m$ edges. If $[z^{2j+1}] g \equiv 0$ for all $j\ge 0$, then
  \begin{eqnarray}
 &{\rm(i)}&  \lambda_1(G_n)  \ge  \frac{m_2}{n-n_1} + \sqrt{\frac{m_2^2}{(n-n_1)^2}+\frac{(m-m_2)^2}{n_1(n-n_1)}}  \nonumber\\[2mm]
 &{\rm(ii)}&\mu_1(G_n)     \ge   \frac{n (m-m_2)}{n_1 (n-n_1)}\nonumber\\[2mm]
 &{\rm(iii)}& q_1(G_n)      \ge  \frac{1}{2} \left( c_1+c_2+\sqrt{\Big(c_1-c_2\Big)^2+\frac{4(m-m_2)^2}{n_1(n-n_1)} }  \right)\nonumber\\[2mm]
 &{\rm(iv)}&  q_2(G_n)      \ge  \frac{1}{2} \left( c_1+c_2-\sqrt{\Big(c_1-c_2\Big)^2+\frac{4(m-m_2)^2}{n_1(n-n_1)} }  \right) \ge q_n(G)\nonumber
  \end{eqnarray} where  $n_1=\lceil\frac{n}{2} \rceil$,  $m_2=m(G_{n-n_1}((\frac{gf}{z})'(\sqrt{z}), f))$,  $c_1=\frac{m-m_2}{n_1}$ and $c_2=\frac{m+3m_2}{n-n_1}$.
  \end{corollary}

 The following result is obtained by Lemma \ref{check} and Corollary  \ref{cor2}.

 \begin{corollary}
 Let $G_n=G_n(g, f)$ be a Riordan graph of the checkerboard type with $n$ vertices and $m$ edges, and let $n_1=\lceil\frac{n}{2} \rceil$. Then,

 \begin{itemize}
 \item[{\rm(i)}] $\lambda_1(G_n)=-\lambda_n(G_n)  \ge  \frac{m}{\sqrt{n_1 (n-n_1)}}.$
 \item[{\rm(ii)}] $\mu_1(G_n)=q_1(G_n)     \ge   \frac{n m}{n_1 (n-n_1)}.$
\item[{\rm(iii)}] $q_n(G_n) = 0$.
 \end{itemize}
  \end{corollary}

 \begin{corollary}
 Let $G_n=G_n(g, f)$ be a Riordan graph with $n$ vertices and $m$ edges such that $[z^{2j}] g=[z^{2j}] f\equiv 0$ for all $j\ge 0$. Then $G_n$ is disconnected with components $H_1\cong G_{n_1}(g'(\sqrt{z}), f)$ and  $H_2\cong G_{n-n_1}((\frac{gf}{z})'(\sqrt{z}), f)$, where $n_1=\lceil\frac{n}{2} \rceil$. Moreover, letting $m_1=m(H_1)$, we have

 \begin{itemize}
\item[{\rm(i)}] $\lambda_1(G_n)=\max\{\lambda_1(H_1), \lambda_1(H_2) \}\ge \max\left\{\frac{2m_1}{n_1}, \frac{2(m-m_1)}{n-n_1}
\right\}$.
\item[{\rm(ii)}] $q_1(G_n)=\max\left\{q_1(H_1), q_1(H_2) \right\}\ge \max\left\{\frac{4m_1}{n_1}, \frac{4(m-m_1)}{n-n_1} \right\}$.
  \end{itemize}
  \end{corollary}

 The following result is obtained by Lemma~\ref{lem12} and Theorem~\ref{co1}.

\begin{corollary}\label{cor-Bell-type}
 Let $G_n=G_n(g(z), zg(z))$ be a Riordan graph of the Bell type with  $n$ vertices and $m$ edges. Then,
  \begin{eqnarray}
 &{\rm(i)}&  \lambda_1(G_n)  \ge  \frac{m_1}{n_1}+ \sqrt{ \frac{m_1^2}{n_1^2}+\frac{(m-m_1)^2}{n_1(n-n_1)}}  \nonumber\\[2mm]
 &{\rm(ii)}&\mu_1(G_n)     \ge   \frac{n (m-m_1)}{n_1 (n-n_1)} \nonumber\\[2mm]
 &{\rm(iii)}& q_1(G_n)      \ge  \frac{1}{2} \left( c_1+c_2+\sqrt{\Big(c_1-c_2\Big)^2+\frac{4(m-m_1)^2}{n_1(n-n_1)} }  \right)\nonumber \\[2mm]
 &{\rm(iv)}&  q_2(G_n)      \ge  \frac{1}{2} \left( c_1+c_2-\sqrt{\Big(c_1-c_2\Big)^2+\frac{4(m-m_1)^2}{n_1(n-n_1)} }  \right) \ge q_n(G)\nonumber
  \end{eqnarray}
  where $n_1=\lceil\frac{n}{2} \rceil$, $c_1=\frac{m+3m_1}{n_1}$, $c_2=\frac{m-m_1}{n-n_1}$ and
$$m_1=\left\{\begin{array}{ll} m(G_{n_1}(g'(\sqrt{z}), zg(z))) & \mbox{ in the general case;} \\
m(G_{n_1}(g(z), zg(z))) & \mbox{ if $G_n$ is io-decomposable.}\end{array}\right.$$
 \end{corollary}

 In the following proposition, we use Lemmas~\ref{PC0} and~\ref{PC}, and Corollary~\ref{cor-Bell-type} to present lower bounds on some eigenvalues and (signless) Laplacian eigenvalues of $PG_{2^k}$, $PG_{2^k+1}$, $CG_{2^k}$ and $CG_{2^k+1}$.

 \begin{proposition}
 For $k\ge 2$, we have
\begin{itemize}
\item[{\rm(i)}] $\lambda_1(PG_{2^k})   \ge \displaystyle{\frac{3^{k-1}+\sqrt{ (3^{k-1}-2^{k-1})^2+4(3^{k-1}-2^{k-2})^2}}{2^{k-1}}-1}
$.
\item[{\rm(ii)}] $  q_1(PG_{2^k})       \ge \displaystyle{\frac{1}{2}\left( \frac{3^{k-1}}{2^{k-4}}+\frac{\sqrt{(3^{k-1}-2^{k-1})^2+(3^{k-1}-2^{k-2})^2}}{2^{k-3}}
\right)-3}$.
\item[{\rm(iii)}]$\lambda_1(PG_{2^k+1}) \ge \displaystyle{\frac{3^{k-1} \left(1+\sqrt{17+2^{3-k}}
\right)}{2^{k-1}+1}}$.
\item[{\rm(iv)}] $ q_1(PG_{2^k+1})      \ge \displaystyle{\frac{3^{k-1}\,\left( 2^{k+1}+1+\sqrt{2^{2k+1}+1}
\right)}{2^{k-1}\,(2^{k-1}+1)}}$.
\item[{\rm(v)}] $ \lambda_1(CG_{2^k})  \ge  \displaystyle{\frac{5\,\cdot\,
3^{2k-2}-3^{k-1}}{2^k}}$.
\item[{\rm(vi)}] $ \mu_1(CG_{2^k})      \ge
\displaystyle{\frac{3^{k-1}}{2^{k-2}}}$.
\item[{\rm(vii)}] $ q_1(CG_{2^k})        \ge
\displaystyle{\frac{2\,\cdot\,3^{k-1}-1+\sqrt{2\,\cdot\,3^{k-1}(3^{k-1}-1)+1}}{2^{k-1}}}$.
\item[{\rm(viii)}] $ q_2(CG_{2^k})        \ge \displaystyle{\frac{2\,\cdot\,3^{k-1}-1-\sqrt{2\,\cdot\,3^{k-1}(3^{k-1}-1)+1}}{2^{k-1}}\ge
q_n(CG_{2^k})}$.
\item[{\rm(ix)}] $\lambda_1(CG_{2^k+1}) \ge  \displaystyle{\frac{b+\sqrt{b^2+(1+2^{1-k})(3^{k-1}+2^{k-1})^2}}{2^{k-1}+1}}$, where
$b=\frac{3^{k-1}+2^k-1}{2}$.
\item[{\rm(x)}] $\mu_1(CG_{2^k+1})     \ge
\displaystyle{\frac{(2^k+1)(3^{k-1}+2^{k-1})}{2^{2k-2}+2^{k-1}}}$.
 \end{itemize}
 \end{proposition}

The following is a result on the spectral radius and the (signless)
Laplacian spectral radius of an io-decomposable Riordan graph
$G_n=G_n(g, zg)$.

\begin{theorem}\label{thm3}
Let $G_n=G_n(g, zg)$ be io-decomposable with $n$ vertices and $m$ edges. Then
\begin{itemize}
\item[{\rm(i)}]  $m\le \frac{n^2}{2}\left( 1-\frac{1}{\lceil\log_2 n\rceil+1}  \right).$
\item[{\rm(ii)}]  $\lceil\log_2 n\rceil \le \lambda_1(G_n)\le n\,\left( 1-\frac{1}{\lceil\log_2 n\rceil+1}  \right).$
\item[{\rm(iii)}]  $\frac{\lambda_1(G_n)}{|\lambda_n(G_n)|}\le \lceil\log_2
n\rceil$.
\item[{\rm(iv)}]  $ \lceil\log_2 n\rceil+1 \le\frac{\lceil\log_2 n\rceil+1}{\lceil\log_2 n\rceil}\lambda_1(G_n)\le
\mu_1(G_n)$.
\item[{\rm(v)}]   $\lambda_1(G_n)-\lambda_1(G_{\lceil \frac{n}{2}\rceil})\le
-\lambda_n(G_n)$.
\item[{\rm(vi)}]   $q_1(G_n) \ge 2\lceil\log_2 n\rceil$.
\end{itemize}
\end{theorem}

\begin{pf}
Considering $X^T=\big[\underbrace{1/n,\, 1/n,\, \ldots,\,
1/n}_n\big]$ in Lemma~\ref{lm1}, and applying Lemma~\ref{lm2}, we
obtain
$$\frac{2}{n^2} \sum_{ij\in E(G_n)}1=\frac{2m}{n^2}\le 1-\frac{1}{\lceil\log_2 n\rceil+1}~\Rightarrow~m\le \frac{n^2}{2}\left(1-\frac{1}{\lceil\log_2 n\rceil+1}  \right).$$ Hence, the proof of (i) is complete. We obtain (ii), (iii) and (iv) directly from Lemmas~\ref{lm2} and~\ref{lm3111}.

Let ${\cal A}_n$ be the adjacency matrix of
$G_n$. By Lemma \ref{lem12} and the definition of an io-decomposable
Riordan graph, we have $P^T {\cal A}_n P=\left[
                                                                                           \begin{array}{cc}
                                                                                             {\cal A}_{\lceil \frac{n}{2}\rceil} & B\\[2mm]
                                                                                             B^T & O \\
                                                                                           \end{array}
                                                                                         \right]$.
From this, after applying Lemma~\ref{par}, we obtain $\lambda_1({\cal
A}_n)+\lambda_n({\cal A}_n)\le \lambda_1( {\cal A}_{\lceil
\frac{n}{2}\rceil})$, which gives the desired result in (v). Finally, (vi) follows from (ii) and (\ref{av2}).
\end{pf}

In what follows, we give an upper bound for the absolute
value of the smallest eigenvalue of an io-decomposable Riordan graph
$G_n$.
\begin{theorem}
Let $G_n=G_n(g, zg)$ be an io-decomposable Riordan graph with the adjacency
matrix ${\cal A}_n:=P {\footnotesize \left[
                     \begin{array}{cc}
                         {\cal A}_{\left\lceil \frac{n}{2} \right\rceil}& B\\[3mm]
                       B^T& $O$\\
                     \end{array}
                   \right]} P^T$, and let $\lambda_i(G_n)$ and  $\sigma_{\max}(B)$ be, respectively, the $i$th largest eigenvalue of $G_n$ and the largest singular value of $B$. Also, let $\eta:=\min\limits_{i=1, \ldots, \lceil \frac{n}{2}\rceil} |\lambda_i(G_{\lceil \frac{n}{2}\rceil})|$.  Then
$$|\lambda_n(G_n)| \le \frac{2\sigma_{\max}^2 (B)}{\eta+\sqrt{\eta^2+4\sigma_{\max}(B)}}.$$
\end{theorem}

\begin{pf} By  Theorem~\ref{thm3} (v) and
Lemma ~\ref{lem61}, we obtain $$|\lambda_n(G_n)|\le
\max_{\substack{i=1,\ldots,\lceil \frac{n}{2}\rceil{}\\j=1+\lceil
\frac{n}{2}\rceil,\ldots,n}}\left\{|\lambda_i(G_n)-\lambda_i(G_{\lceil
\frac{n}{2}\rceil})|\, ,\, |\lambda_j(G_n)|\right\}\le
\frac{2\sigma_{\max}^2 (B)}{\eta+\sqrt{\eta^2+4\sigma_{\max}(B)}}.$$
\end{pf}

The following result is related to the structure of the adjacency
matrix of an io-decomposable Riordan graph $G_n(g,zg)$, and it
follows from the definition of an io-decomposable Riordan graph of
the Bell type and Lemma \ref{lem12}.

 \begin{lemma}\label{lm3}
 Let ${\cal A}_n:= P {\footnotesize \left[
                     \begin{array}{cc}
                         {\cal A}_{\left\lceil \frac{n}{2} \right\rceil}& B\\[3mm]
                       B^T& $O$\\
                     \end{array}
                   \right]} P^T$ be the adjacency matrix of an io-decomposable Riordan graph $G_n(g,zg)$, and let $b_{i,j}$  be the $(i,j)$th entry of $B^T$. If $n$ is odd, then $P^{T}{\cal A}_n P$ is given by

\begin{footnotesize}
 \begin{eqnarray*}{\displaystyle
 \left[
           \begin{array}{cccc|ccc}
             0 & r_{1,2}  & \cdots & r_{1,k} & b_{1,1} & \cdots & b_{k-1,1}\\
             r_{2,1} & 0  & \ddots & \vdots &  r_{2,1} & \ddots & \vdots \\
             \vdots &\ddots  &\ddots & r_{k-1,k}  &  \vdots & \ddots & b_{k-1,k-1} \\
                            r_{k,1} &\cdots   & r_{k,k-1}& 0 & r_{k,1}  & \cdots & r_{k,k-1} \\
               \hline
             b_{1,1} & r_{1,2}  & \cdots & r_{1,k} & 0 & \cdots & 0 \\
             \vdots  & \ddots & \ddots & \vdots & \vdots & \ddots & \vdots \\
             b_{k-1,1} &\cdots & b_{k-1,k-1} & r_{k-1,k} &0 & \cdots & 0 \\
           \end{array}
         \right]
}
  \end{eqnarray*} \end{footnotesize} where $k=\frac{n+1}{2}$. Otherwise, if $n$ is even and $j=\frac{n}{2}$, then $P^{T}{\cal A}_n P$ is

\begin{footnotesize}\begin{eqnarray*}
{\displaystyle
 \left[
           \begin{array}{cccc|cccc}
             0 & r_{1,2}  & \cdots & r_{1,j} & b_{1,1} & \cdots & \cdots& b_{j,1}\\
             r_{2,1} & 0 ~~ &\ddots & \vdots &  r_{2,1}  & \ddots & \vdots & \vdots \\
             \vdots & \ddots &\ddots & r_{j-1,j} &   \vdots &  \ddots &\ddots  & \vdots\\
             r_{j,1} & \cdots  & r_{j,j-1}& 0 &r_{j,1} & \cdots & r_{j,j-1}& b_{j,j} \\
               \hline
             b_{1,1} & r_{1,2} & \cdots & r_{1,j} & 0 & \cdots & \cdots&0 \\
             \vdots & \ddots & \ddots & \vdots & \vdots & \ddots & \vdots & \vdots \\
             \vdots & \vdots & \ddots & r_{j-1,j} & \vdots  & \ddots & \ddots &\vdots\\
             b_{j,1}  & \cdots &\cdots &  b_{j,j} & 0  & \cdots & \cdots&0 \\
           \end{array}
         \right].
}
 \end{eqnarray*}
 \end{footnotesize}
 \end{lemma}

 In the following, we
give a lower bound on the spectral radius of any io-decomposable
Riordan graph $G_n(g,zg)$. We shall mean
$$\left\{gf^k\right\}_n(1):=\sum_{j=k}^n[z^j]\left\{gf^k\right\}$$ by
the substitution of $z=1$ in the Taylor expansion in $z$ of
$\left\{gf^k\right\}$ up to degree $n$ modulo 2.

\begin{theorem}
Let $G_n=G_n(g,zg)$ be an io-decomposable Riordan graph. Then,
$$\lambda_1(G_n)\ge \frac{1+\sqrt{2}}{n_1}\,\sum_{j=0}^{n_1-2} \{
g^{\prime}(\sqrt{z}) \cdot f^j\}_{n_1-2}(1).$$
\end{theorem}
\begin{pf} By Corollary~\ref{cor-Bell-type}, we have
\begin{equation}\label{ahm3}
\lambda_1(G_n) \ge
\frac{m_1}{n_1}+\sqrt{\frac{m_1^2}{n_1^2}+\frac{(m-m_1)^2}{n_1(n-n_2)}}.
\end{equation} Let ${\cal A}_n=P {\footnotesize \left[
                     \begin{array}{cc}
                         {\cal A}_{\left\lceil \frac{n}{2} \right\rceil}& B\\[3mm]
                       B^T& $O$\\
                     \end{array}
                   \right]} P^T$ be the adjacency matrix of $G_n$. Then, by Lemma \ref{lm3},
                   we have  $$\sigma(B)>m\left({\cal A}_{\left\lceil \frac{n}{2} \right\rceil}\right)=m_1\ge \sqrt{\frac{n-n_1}{n_1}}\,m_1~\Rightarrow~\sigma(B)^2>\frac{n-n_1}{n_1} m_1^2\Rightarrow\frac{(m-m_1)^2}{n_1(n-n_1)}>\frac{m_1^2}{n_1^2}.$$
From (\ref{ahm3}), the above and Lemma \ref{lem12}, the required result is obtained.
\end{pf}

 The following
two results give lower bounds on the largest Laplacian eigenvalue
$\mu_1(G_n)$ of a Riordan graph $G_n$.
 \begin{theorem} \label{lap1}
 Let $G_n=G_n(g, f)$ be a Riordan graph with $n$ vertices and $m$ edges, and let  $n_1=\lceil\frac{n}{2} \rceil$ and $n_2=n-n_1$. Then
 \begin{eqnarray}
 \mu_1(G_n)\ge\frac{n}{n_1 n_2}\left( \sum_{j=0}^{n_2-1}\left\{z(gf)^{\prime}(\sqrt{z})\cdot f^j\right\}_{n_1-1}(1)+\sum_{j=0}^{n_1-1}\left\{(zg)^{\prime}(\sqrt{z})\cdot f^j\right\}_{n_2-1}(1)\right).\nonumber
  \end{eqnarray}
  \end{theorem}
\begin{pf}
%
Let $L=D-\mathcal{A}$ be the Laplacian matrix of $G_n$. Using the
same technique as in the proof of Theorem~\ref{thm9}, and by
Lemma~\ref{lem31} applied to $PLP^T$, we obtain
 $$B_2:=\left[
           \begin{array}{cc}
            \frac{\sigma(B)}{n_1} & -\frac{\sigma(B)}{n_1} \\
                             &                  \\
             -\frac{\sigma(B)}{n_2} &  \frac{\sigma(B)}{n_2} \\
           \end{array}
         \right],$$ where $\sigma(B)$ denotes the number of 1s in the matrix B.
 It is easy to see that the eigenvalues of $B_2$ are given by
 $$\lambda_1(B_2)=\frac{n\sigma(B)}{n_1 n_2}  ~~~\mbox{and}~~~ \lambda_2(B_2)=0.$$ By Lemma \ref{lem31}, we have
\begin{equation}
\mu_1(G)\ge \lambda_1(B_2)=\frac{n\sigma(B)}{n_1 n_2}. \label{mu1}
\end{equation}

Now, by (iv) in the Riordan Graph Decomposition Theorem (Theorem~\ref{e:th}), we have
\begin{eqnarray}
\sigma(B)&=&\sigma\left( {\cal B}(z\cdot(gf)^{\prime }(\sqrt{z}),f(z))_{n_1\times n_2}\right)+\sigma\left({\cal B}((zg)^{\prime}(\sqrt{z}),f(z))_{n_2 \times n_1}^{T}\right)\nonumber\\[2mm]
         &=&\sigma\left( {\cal B}(z\cdot(gf)^{\prime }(\sqrt{z}),f(z))_{n_1\times n_2}\right)+\sigma\left({\cal B}((zg)^{\prime}(\sqrt{z}),f(z))_{n_2 \times n_1}\right)\nonumber\\[2mm]
         &=&\sum_{j=0}^{n_2-1}\left\{z(gf)^{\prime}(\sqrt{z})\cdot f^j\right\}_{n_1-1}(1)+\sum_{j=0}^{n_1-1}\left\{(zg)^{\prime}(\sqrt{z})\cdot f^j\right\}_{n_2-1}(1).\nonumber
\end{eqnarray}
 Using this in (\ref{mu1}), the required result is obtained.
\end{pf}

The following result is an immediate corollary of Theorem~\ref{lap1}.

 \begin{corollary}\label{lap2}
 Let $G_n=G_n(g, f)$ be a proper Riordan graph with $n$ vertices and $m$ edges. Then
 \begin{eqnarray}
 \mu_1(G_n)\ge \frac{n}{n_1 n_2}\left( \left\{z(gf)^{\prime}(\sqrt{z})\right\}_{n_1-1}(1)+\left\{(zg)^{\prime}(\sqrt{z})\right\}_{n_2-1}(1)+2\left\lfloor \frac{n}{2}\right\rfloor-3 \right).\nonumber
 \end{eqnarray}
 \end{corollary}

 For a graph $G$, let $a(G)$ be the {\em algebraic connectivity} of $G$. The following result gives a relation between the algebraic connectivity of $\langle V_o\rangle$ and the median Laplacian eigenvalues
 $\mu_{\lceil \frac{n}{2}\rceil}(G_n)$ and $\mu_{\lceil \frac{n}{2}\rceil+1}(G_n)$ of  an o-decomposable Riordan graph $G_n$.

\begin{theorem} \label{tm21}
Let $G_n=G_n(g ,f)$ be an o-decomposable Riordan graph and let
$H\cong \langle V_o\rangle$. Then
\begin{align}
   \mu_{\left\lceil \frac{n}{2}\right\rceil+1}(G_n)&\le \left\lceil \frac{n}{2}\right\rceil   \label{equ111}\\[2mm]
   \mu_{\left\lceil \frac{n}{2}\right\rceil}(G_n)  &\le\left\{
                                               \begin{array}{ll}
                                                  \lfloor \frac{n}{2}\rfloor+a(H) & \hbox{if $n$ is even or $n$ is odd with $a(H)> 1$;} \\[3mm]
                                                 \lceil \frac{n}{2}\rceil & \hbox{if $n$ is odd with $a(H)\le 1$.}
                                               \end{array}
                                             \right.
    \label{equ222}
\end{align}
\end{theorem}

\begin{pf}
One can see that $K_{\left\lfloor \frac{n}{2} \right\rfloor}$ is a subgraph of $\overline{G}_n$. Further, it is well known \cite{cv} that the Laplacian spectrum of $K_{\left\lfloor \frac{n}{2} \right\rfloor}$ is as follows:
$$Lspec\left(K_{\left\lfloor \frac{n}{2} \right\rfloor}\right)=\Big\{\underbrace{\left\lfloor \frac{n}{2} \right\rfloor,\ldots, \left\lfloor \frac{n}{2} \right\rfloor}_{\left\lfloor \frac{n}{2} \right\rfloor-1},0   \Big\}.$$
From the above, we obtain $$\mu_{\left\lfloor
\frac{n}{2}\right\rfloor-1}(\overline{G}_n)\ge \mu_{\left\lfloor
\frac{n}{2}\right\rfloor-1}\left(K_{\left\lfloor \frac{n}{2} \right\rfloor}\right)=\left\lfloor
\frac{n}{2} \right\rfloor.$$ From this and Lemma~\ref{lem51}, we obtain
$$n-\mu_{\left\lceil \frac{n}{2}\right\rceil+1}(G_n)\ge \left\lfloor \frac{n}{2} \right\rfloor,$$ which completes the proof of (\ref{equ111}).

Now, we consider the subgraph $\overline{H}\cup K_{\left\lfloor
\frac{n}{2} \right\rfloor}$ of $\overline{G}_n$,   where $H \cong \langle
V_o\rangle$. By Lemma \ref{lem51}, we have
$\mu_i(\overline{H})=\left\lceil \frac{n}{2}\right\rceil-\mu_{\left\lceil
\frac{n}{2} \right\rceil-i}(H)$ for $1\le i \le \left\lceil
\frac{n}{2}\right\rceil-1$. From this and Lemma~\ref{lem51}, we have
\begin{equation}
\mu_{\left\lfloor \frac{n}{2} \right\rfloor}(\overline{G}_n)\ge
\min\Big\{\left\lfloor \frac{n}{2} \right\rfloor,
\mu_1(\overline{H})\Big\}=\min\Big\{\left\lfloor \frac{n}{2}
\right\rfloor,\left\lceil  \frac{n}{2}\right\rceil- a(H)\Big\}.  \label{equ3}
\end{equation} If $n$ is even, from (\ref{equ3}) and the fact  that $\mu_1(\overline{H})\le |V(\overline{H})|=\frac{n}{2}$, we obtain
$$\mu_{\frac{n}{2}}(\overline{G}_n)\ge \mu_1(\overline{H})~\Rightarrow~n-\mu_{\frac{n}{2}}(G_n)\ge \frac{n}{2}-\mu_{\frac{n}{2}-1}(H)
=\frac{n}{2}-a(H)~\Rightarrow~\mu_{\frac{n}{2}}(G_n)\le
\frac{n}{2}+a(H),$$ which gives the result in (\ref{equ222}) for
even $n$.  Otherwise, $n$ is odd. First, we assume that
$\mu_1(\overline{H})\ge \lfloor \frac{n}{2} \rfloor$, that is,
$a(H)\le 1$. From this and (\ref{equ3}), we
have$$\mu_{\frac{n-1}{2}}(\overline{G}_n)\ge
\frac{n-1}{2}~\Rightarrow~\mu_{\frac{n+1}{2}}(G_n)\le
\frac{n+1}{2}.$$ Next, we assume that $\mu_1(\overline{H})< \lfloor
\frac{n}{2} \rfloor$, that is, $a(H)>1$. From this and (\ref{equ3}),
we have $$\mu_{\frac{n-1}{2}}(\overline{G}_n)\ge
\mu_1(\overline{H})~\Rightarrow~\mu_{\frac{n+1}{2}}(G_n)\le
\frac{n-1}{2}+a(H).$$ This completes the proof.\end{pf}

By similar arguments as those in Theorem~\ref{tm21}, one can
easily prove the following theorem on e-decomposable Riordan graphs.
\begin{theorem} \label{tm31}
Let $G_n=G_n(g ,f)$ be an e-decomposable Riordan graph and let
$H\cong \langle V_e\rangle$. Then
\begin{align*}
   \mu_{\left\lfloor \frac{n}{2}\right\rfloor+1}(G_n)\le \left\lfloor \frac{n}{2}\right\rfloor~{\rm and}~\mu_{\left\lfloor \frac{n}{2}\right\rfloor}(G_n)  \le \left\lfloor \frac{n}{2}\right\rfloor+a(H).
\end{align*}
\end{theorem}

\section{Vertex degrees in Riordan graphs}\label{vertexdegrees-sec}

From now on, we assume that $p=\lfloor\log_2(n-1)\rfloor$ and
$\{g\}_{-1}(1)=0$. In the following, we obtain the degree of the
vertex $2^p+1$ in an io-decomposable Riordan graph $G_n(g ,zg)$.

\begin{lemma}\label{imp}
Let $G_n=G_n(g ,zg)$ be an io-decomposable Riordan graph. Then
\begin{equation}\label{ahm4}
d_{G_n} (2^p+1)=2^p+\{g\}_{n-2^p-2}(1).
\end{equation}  In particular, the vertex $2^p+1$ is universal in $G_n$ if and only if $n=2^p+1$ or $n\neq 2^p+1$ and $[z^{2s}] g\equiv 1$ for $0\le s \le \lfloor \frac{n-2^p}{2}\rfloor-1$.
\end{lemma}
\begin{pf}
 Suppose that $H=\langle \{2^p+1,\ldots,n\}\rangle$. Then
\begin{align*}
d_{G_n} (2^p+1)&= d_{G_{2^p+1}}(2^p+1)+d_H(2^p+1)\\[2mm]
               &= d_{G_{2^p+1}}(2^p+1)+d_{G_{n-2^p}}(1) ~~~~~~\mbox{(by Theorem 3.6(i) in \cite{CJKM}) }\\[2mm]
               &= 2^p+\{g\}_{n-2^p-2}(1) ~~~~~~~~\mbox{(by Theorems 4.13 and 3.1 in \cite{CJKM})}.
\end{align*} This completes the proof for (\ref{ahm4}). The second part of the proof follows from (\ref{ahm4}) and Lemma 4.4 in \cite{CJKM}.
\end{pf}

For any io-decomposable Riordan graph $G_n=G_n(g ,zg)$, we have
$[z^0]g=[z^1]g\equiv 1$. This observation, along with
Lemma~\ref{ahm4}, gives the following result.
\begin{corollary} \label{cor5}
Let $G_n=G_n(g ,zg)$ be an io-decomposable Riordan graph. Then the
vertex $2^p+1$ is universal in $G_n$ for $2^p+1\le n\le 2^p+3$.
\end{corollary}

In the following two lemmas, some comparisons between the degree of
the vertex $2^p+1$ and degrees of other vertices in an
io-decomposable Riordan graph are given.
\begin{lemma}\label{ahm7}
Let $G_n=G_n(g ,zg)$ be an io-decomposable Riordan graph. Then
\begin{equation}
d_{G_n}(2^p+1)\ge d_{G_n}(1).  \label{equ4}
\end{equation}
Moreover, for $n\le 1+2^{p-1}+2^p$, we have
\begin{equation} \label{equ5}
d_{G_n}(2^p+1)\ge d_{G_n}(2^{p-1}+1).
\end{equation}
\end{lemma}
\begin{pf}
From Lemma \ref{imp}, we have
\begin{eqnarray*}
d_{G_n}(2^p+1)=2^p+\{g\}_{n-2^p-2}(1)\ge\sum_{i=n-2^p-1}^{n-2} [z^i]
g+\{g\}_{n-2^p-2}(1)=\{g\}_{n-2}(1)=d_{G_n}(1).
\end{eqnarray*} This completes the proof of (\ref{equ4}).

Now, suppose that  $n\le 1+2^{p-1}+2^p$. Then
\begin{equation}\label{equ6}
\{g\}_{2^{p-1}-1}(1)-\{g\}_{n-2^p-2}(1)=\sum_{i=n-2^p-1}^{2^{p-1}-1}
[z^i] g\le 2^p+2^{p-1}+1-n.
\end{equation} We define the subgraphs $H_1$ and $H_2$ of $G_n$ as follows:
$$H_1\cong \left\{
    \begin{array}{ll}
      \langle \{2^{p-1}+1, 2^p+2, \ldots, n\}\rangle & \hbox{if $n>2^p+1$} \\[3mm]
       \langle \{2^{p-1}+1\}\rangle & \hbox{if $n=2^p+1$}
    \end{array}
  \right.$$  and  $H_2\cong \langle \{2^{p-1}+1,\ldots, 2^p+1\}\rangle$.
 By Lemma \ref{imp}, we have
\begin{align*}
d_{G_n}(2^p+1)&=2^p+\{g\}_{n-2^p-2}(1)\\
              &\ge n-2^{p-1}-1+\{g\}_{2^{p-1}-1}(1)~~~~~~\mbox{(by ~(\ref{equ6}))}\\
              &=2^{p-1}+n-2^p-1+\{g\}_{2^{p-1}-1}(1)\\
              &\ge2^{p-1}+\{g\}_{2^{p-1}-1}(1)+d_{H_1}(2^{p-1}+1)\\
              &=d_{G_{2^{p-1}+1}}(2^{p-1}+1)+d_{G_{2^{p-1}+1}}(1)+d_{H_1}(2^{p-1}+1)\\
              &=d_{G_{2^{p-1}+1}}(2^{p-1}+1)+d_{H_2}(2^{p-1}+1)+d_{H_1}(2^{p-1}+1)\\
              &=d_{G_n}(2^{p-1}+1),
\end{align*} where the next to last equality is given by Theorem 3.6 in \cite{CJKM}, and the last inequality is obtained using the fact that $d_{H_1}(2^{p-1}+1)\le n-2^p-1$. This completes the proof of (\ref{equ5}).
\end{pf}

\begin{lemma}
Let $G_n=G_n(g ,zg)$ be an io-decomposable Riordan graph. Then, for
$i\in V_e$,
$$d_{G_n}(2^p+1)\ge d_{G_n}(i).$$
\end{lemma}

\begin{pf}
We have $n\le 2^{p+1}$, that is, $ \left\lceil
\frac{n}{2}\right\rceil \le 2^p$. From this and (\ref{ahm4}), for
$i\in V_e$, we obtain
$$d_{G_n}(2^p+1)\ge 2^p \ge \left\lceil \frac{n}{2}\right\rceil\ge d_{G_n}(i).$$
\end{pf}

\begin{corollary}
Let $G_n=G_n(g,zg)$ be an io-decomposable Riordan graph. Then
$$2^p+\left\lceil \log_2 (n-2^p)\right\rceil\le d_{G_n}(2^p+1)\le n-1.$$
\end{corollary}

\begin{pf}
By Lemma~\ref{imp}, we have
\begin{equation}\label{equ7}
d_{G_n} (2^p+1)=2^p+\{g\}_{n-2^p-2}(1).
\end{equation} Since $G_n$ is an io-decomposable Riordan graph, by Lemma~\ref{e:lem3}, we have $g_{2^p-1}\equiv 1$ for $p\ge 0$.
From this, (\ref{equ7}) and Lemma \ref{lm2}, we obtain $$d_{G_n}
(2^p+1)=2^p+\{g\}_{n-2^p-2}(1)\ge 2^p+\lceil \log_2 (n-2^p)\rceil.$$
This completes our proof.
\end{pf}

\begin{lemma}[\cite{CJKM}]\label{moj3}
Let $G_{n}(g,zg)$ be an io-decomposable Riordan graph. Then it is
$(\lceil\log_2n\rceil+1)$-partite with partitions
$V_1,V_2,\ldots,V_{\lceil\log_2n\rceil+1}$ such that

$$V_j=\left\{2^{j-1}+1+i2^j\,|\,0\le i\le\left\lfloor
{n-1-2^{j-1}\over 2^j}\right\rfloor\right\}~~~ \text{if $1\le
j\le\lceil\log_2n\rceil$},$$
$$V_{\lceil\log_2n\rceil+1}=\{1\}.~~~~~~~~~~~~~~~~~~~~~~~~~~~~~~~~~~~~~~~~~~~~~~~~~~~~~~~~~~~~~~~~~$$
\vskip1pc
\end{lemma}

It follows from Lemma~\ref{ahm7} that if the vertex $1$ is universal
in an io-decomposable Riordan graph $G_n(g ,zg)$, then the vertex
$2^p+1$ is also universal in that graph. In what  follows, we show
that if the vertex $2^{p-1}+1$ is universal, then the vertex $2^p+1$
is also universal.

\begin{theorem}
Let $G_n=G_n(g,zg)$ be an io-decomposable Riordan graph. If
$2^{p-1}+1$ is universal then the vertex $2^p+1$ is also universal.
\end{theorem}

\begin{pf}
By Lemma~\ref{ahm7}, it suffices to prove that if $2^{p-1}+1$ is a
universal vertex in $G_n$, then $n\le 1+2^{p-1}+2^p$. By Lemma
\ref{moj3}, we obtain that the vertex $2^{\left\lceil \log_2
n\right\rceil-2}+1\in V_{\left\lceil\log_2 n\right\rceil-1}$.

We claim that $p=\lceil \log_2 n\rceil-1$. Indeed, by definition of
$p$, we have $2^p+1\le n \le 2^{p+1}$. We may assume that
$n=2^p+\gamma$, where $1\le \gamma \le2^p$. So,
$$p=\left\lceil \log_2 (2^p+1)\right\rceil-1\le\lceil \log_2 n\rceil-1=\left\lceil \log_2 (2^p+\gamma)\right\rceil-1\le \left\lceil \log_2 2^{p+1}\right\rceil-1=p,$$
which completes the proof of our claim.

Thus, $2^{p-1}+1\in V_{p}$. Again, by Lemma~\ref{moj3}, we have
$|V_p|=\left\lfloor \frac{n-1-2^{p-1}}{2^p}\right\rfloor +1$. Since
the vertex $2^{p-1}+1$ is universal, we must have $|V_p|=1$ and,
consequently, $\left\lfloor \frac{n-1-2^{p-1}}{2^p}\right\rfloor=0$,
so $\frac{n-1-2^{p-1}}{2^p}<1$ and $n< 1+2^{p-1}+2^p$. This
completes the proof.
\end{pf}

Based on the facts above, we state the following conjecture.
\begin{conjecture}
 The vertex $2^p+1$ has the maximum degree in any io-decomposable Riordan graph
 $G_n=G_n(g,zg)$.
\end{conjecture}

The following result is obtained from the
facts above.
\begin{corollary}\label{cor-4-9}
Let $G_n=G_n(g,zg)$ be an io-decomposable Riordan graph. Then $G_n$
has at most three universal vertices which must come from the set
$\{1,2^{p-1}+1,2^{p}+1\}$. In particular, if $n\ge 1+2^{p-1}+2^p$,
then $G_n$ has at most two universal vertices that must come from
the set $\{1,2^p+1\}$.
\end{corollary}

The following result on Pascal graphs is obtained in \cite{CKM}.

 \begin{lemma}[\cite{CKM}]\label{dg7}
For a Pascal graph $PG_n$ of order $n$, we have
\begin{itemize}
 \item[{\rm(i)}] If $n=2^p+1$, then the universal  vertices are $1,\frac{n+1}{2},n$.
 \item[{\rm(ii)}] If $n\neq 2^p+1$, then the universal  vertices are $1,2^{p}+1$.
 \end{itemize}
 \end{lemma}

By Corollary~\ref{cor-4-9}, for any io-decomposable Riordan graph of the
Bell type with three universal vertices, the vertex $1$ must be a
universal vertex. The latter implies that $g(z)=\frac{1}{1-z}$ and
$G_n\cong PG_n$, and thus, by Lemma~\ref{dg7}, we have the following
result.

\begin{theorem} Let $G_n=G_n(g,zg)$ be an io-decomposable Riordan graph
with three universal vertices. Then $G_n\cong PG_n$ and $n=2^p+1$.
\end{theorem}

\begin{corollary}
Let $G_n=G_n(g,zg)$ be an io-decomposable Riordan graph with exactly
one universal vertex. Then the universal vertex must be $2^p+1$.
\end{corollary}

\begin{corollary}
Let $G_n=G_n(g,zg)$ be an io-decomposable Riordan graph. Then at
most three Laplacian eigenvalues of $G_n$ are equal to $n$. In
particular, if $n\ge 1+2^{p-1}+2^p$, then at most two Laplacian
eigenvalues of $G_n$ are equal to $n$.
\end{corollary}

Theorems~\ref{thm-AM1} and~\ref{th41} below give bounds for
$\mu_{1}(G_n)$ and $\mu_{n-1}(G_n)$, respectively, in terms of the
generating function $g(z)$ for any io-decomposable Riordan graph of
the Bell type.

\begin{theorem}\label{thm-AM1}
Let $G_n=G_n(g ,zg)$ be an io-decomposable Riordan graph.  Then
\begin{equation}
\mu_{1}(G_n)\ge 2^p+\{g\}_{n-2^p-2}(1)+1. \label{ahm5}
\end{equation} The equality holds if and only if $n=2^p+1$, or $n\neq 2^p+1$ and $[z^{2s}] g\equiv 1$ for $0\le s \le \lfloor \frac{n-2^p}{2}\rfloor-1$.
\end{theorem}

\begin{pf}
By Lemmas ~\ref{imp} and~\ref{max1}, we obtain
\begin{equation}
\mu_1(G_n)\ge \Delta(G_n)+1\ge d(2^p+1)+1=2^p+\{g\}_{n-2^p-2}(1)+1.
\label{ahm6}
\end{equation}
This completes the proof of (\ref{ahm5}). To prove the rest, we first assume that $n=2^p+1$. In this case, $\Delta(G_n)=n-1$ and
we obtain $\mu_1(G_n)\ge n$. From this and Lemma~\ref{lem51}, we obtain
$\mu_1(G_n)=n$. Otherwise, $n\neq 2^p+1$ and $[z^{2s}] g\equiv 1$ for
$0\le s \le \lfloor \frac{n-2^p}{2}\rfloor-1$. From this, we obtain
$\{g\}_{n-2^p-2}(1)=n-2^p-1$, and thus $2^p+\{g\}_{n-2^p-2}(1)+1=n$.
This gives $\mu_1(G_n)\ge n$ and, again, by Lemma~\ref{lem51}, we obtain
$\mu_1(G_n)=n=2^p+\{g\}_{n-2^p-2}(1)+1$.

Conversely, we assume that $\mu_1(G_n)=2^p+\{g\}_{n-2^p-2}(1)+1$. By (\ref{ahm6}), $$\mu_1(G_n)=\Delta(G_n)+1=
d(2^p+1)+1.$$ But then, by  Lemma~\ref{max1},
$d(2^p+1)=\Delta(G_n)=n-1$, and thus, by Lemma~\ref{imp}, we obtain
$n=2^p+1$, or $n\neq 2^p+1$ and $[z^{2s}] g\equiv 1$ for $0\le s \le
\lfloor \frac{n-2^p}{2}\rfloor-1$. We are done.
\end{pf}

\begin{theorem} \label{th41}
Let $G_n=G_n(g ,zg)$ be an io-decomposable Riordan graph.  Then
\begin{eqnarray*}
\mu_{n-1}(G_n)\leq 1+\sum_{i=1}^{\lfloor\frac{n+1}{2}\rfloor}
g_{2i-1}.
\end{eqnarray*} Moreover, for $n\in\{2^p+1,2^p+2\}$, we have $\mu_{n-1}(G_n)\ge 1$.
\end{theorem}
\begin{pf}
By Lemma \ref{lm3}, we obtain $$d_{G_n}(2)=1+d_{G_{\left\lceil
\frac{n}{2}\right\rceil} }(1)=1+\sum_{i=1}^{\left\lfloor\frac{n+1}{2}\right\rfloor}
g_{2i-1}.$$
 Since $G_n\ncong K_n$, by Lemma~\ref{mb1}, we obtain $$\mu_{n-1}(G_n)\leq \delta(G_n)\le d_{G_n}(2)=1+\sum_{i=1}^{\lfloor\frac{n+1}{2}\rfloor} g_{2i-1},$$
as desired. Now, for $n\in\{2^p+1,2^p+2\}$, by Corollary~\ref{cor5}, we have that $G_n$ has at least one universal vertex. This means that the complement graph $\overline{G_n}$ of the io-decomposable Riordan graph $G_n$ has at least one isolated vertex. Thus, $\mu_1(\overline{G_n})\le
n-1$ and $\mu_{n-1}(G)\ge 1$.
\end{pf}

We end this section by using generating functions to study the {\em
Rayleigh-Ritz quotient} for Riordan graphs. It is well known that
for a graph $G$, the Rayleigh-Ritz quotient under the adjacency
matrix
 $\mathcal{A}(G)$ provides a lower bound on the spectral radius $\lambda_1(G)$ of $G$ as follows.
 \begin{lemma}[\cite{BH}]\label{Ray} Let $G$ be a graph with $n$ vertices and with the adjacency matrix $\mathcal{A}(G)$. For any nonzero vector
 $X=\left[x_1,\,x_2,\ldots,x_n\right]^T \in \mathbb{R}^n$,
 \begin{equation}  \nonumber
 \lambda_1(G)\geq \frac{X^t \mathcal{A}(G) X}{X^t X}.
 \end{equation} The equality holds if and only if $X$ is the eigenvector of $\mathcal{A}(G)$ with the spectral radius $\lambda_1(G)$.
 \end{lemma}

 Let $k=k(z):=k_0+k_1z+\cdots+k_nz^n$ and $\ell=\ell(z):=\ell_0+\ell_1z+\cdots+\ell_nz^n$ be two polynomials of degree $n$ in $R[[z]]$.
 Similarly to the inner product between two vectors,  we define $$k\odot \ell=k(z)\odot \ell(z):=k_0\ell_0+k_1\ell_1+\cdots+k_n\ell_n.$$
The following result gives a lower bound for the spectral radius of a Riordan graph of the Appell type.

 \begin{theorem} \label{th3a}
 Let $G_n=G_n(g,z)$ be a Riordan graph of the Appell type with $n$ vertices. Then for any polynomial $h$ of degree $n-1$, we have
 \begin{equation} \label{moj2}
 \lambda_1(G_n)\geq \frac{2\,\left(h\odot k\right)}{h\odot h},
 \end{equation} where $k=\sum_{i=0}^{n-2} ([z^i]gh) z^i$.
 The equality in (\ref{moj2}) holds  if and only if $h$ is the generating function of the eigenvector of $G_n$ with eigenvalue $\lambda_1(G_n)$.
 \end{theorem}
 \begin{pf} Let $\mathcal{A}_n$ be the adjacency matrix of a Riordan graph $G_n(g,f)$. By Lemma~\ref{Ray} and (\ref{def1}),  for any nonzero vector $X\in R^n$, we have
 \begin{align*}
 \lambda_1(G_n(g,f))\geq\frac{X^T\, {\cal A}_n\, X}{X^T X} &= \frac{X^T \left((zg, f)_n+(zg, f)^T_n\right)X}{X^T X}\\
 &=\frac{X^T(zg, f)_n X+X^T(zg, f)^T_n X}{X^T X}.
 \end{align*}
 Since $\left(X^T(zg, f)_n X\right)^T=X^T(zg, f)^T_n X$ and $X^T (zg, f)_n X\in
 R$, we obtain
 $$X^T (zg, f)_n X=X^T (zg, f)^T_n X~{\rm and}~
 \lambda_1(G_n(g,f))\geq 2\frac{X^T (zg, f)_n X}{X^T X}.$$ Now, we let
 $Y:=(zg, f)_n X$ and $\ell(z)$ and $h(z)$ be, respectively, the generating functions of $Y$ and $X$. Then, for $f=z$, we have
 $$\lambda_1(G_n)\geq 2\frac{X^T Y}{X^T X}=2\frac{h\odot \ell}{h\odot h}=2\frac{h\odot (zgh(f))}{h\odot h}=2\frac{h\odot (zgh)}{h\odot h},$$ where we have applied the FTRM (Lemma~\ref{FTRM-lem}) and the fact that $zgf^j=z^{j+1}g\in Z_2[[z]]$ for any $0\leq j\leq n-2$. This completes the proof of (\ref{moj2}). The second part of the proof follows directly from the equality condition in Lemma~\ref{Ray}. \end{pf}

Next, for the sake of example, we demonstrate that Theorem~\ref{th3a} works for complete graphs.
 \begin{example}
 It is well known that  for a complete graph $K_n$, $\lambda_1(K_n)=n-1$ with eigenvector $X^T=[\underbrace{1, \ldots, 1}_n]$. The generating function of $X$ is
 $h(z)=\sum_{i=0}^{n-1} z^i$. Moreover, by Theorem \ref{th3a}, for $K_n\cong G_n\left(\frac{1}{1-z},z\right)$, we have
 $$h\odot h=\sum_{i=0}^{n-1} 1=n,~~k(z)=\sum_{i=1}^{n-1} iz^{i-1},~~h\odot k=\sum_{i=1}^{n-1}i=\frac{n\,(n-1)}{2},$$ and hence, by Theorem~\ref{th3a}, indeed  $\lambda_1(K_n)=n-1$.
 \end{example}

\section{Eigenvalues of Pascal graphs
and Catalan graphs}\label{PasCat}

The following conjecture shows
the significance of Pascal graphs $PG_n$ and Catalan graphs
$CG_n$, and we devote this section to studying eigenvalues of these graphs.
\begin{conjecture}[\cite{CJKM}]\label{conj} {\rm Let $G_n$ be an
io-decomposable Riordan graph of the Bell type. Then,
\begin{align*}
2={\rm diam}(PG_n)\le {\rm diam}(G_n) \le {\rm diam}(CG_n)
\end{align*}
for $n\geq 4$. Moreover, $PG_n$ is the only graph in the class of
io-decomposable graphs of the Bell type  whose diameter is $2$ for
{\em all} $n\ge4$.}
\end{conjecture}

In this section, we first present some results related to
eigenvalues $-1$ and $0$ of Pascal and Catalan graphs, and next, we give some results about integral Laplacian eigenvalues of these
graphs.

The following result can be directly proved by definition of an
eigenvalue and a corresponding eigenvector of a graph.
\begin{lemma} \label{lem5}
Let $G$ be a graph of order $n$. If $uv\in E(G)$ with
$N_G(u)\backslash  \{v\}=N_G(v)\backslash  \{u\}$, then $-1$ is an
eigenvalue of $G$ with eigenvector $$X=[0,\ldots,0,1\,,0,
\ldots,0,-1\,,0,\ldots,0]^T,$$ where the entries $1$ and $-1$
correspond to the vertices $u$ and $v$, respectively.
\end{lemma}

Recall that the Catalan graph $CG_n$ is defined by
$G_n\left({1-\sqrt{1-4z}\over 2z},{1-\sqrt{1-4z}\over 2}\right)$.
\begin{theorem}\label{thm13}
Let  $n\ge 2$. Then $-1$ is an eigenvalue of a  Catalan graph $CG_n$
with eigenvector $X=[1,-1,0,\ldots,0]^T$.
\end{theorem}
\begin{pf}
We know that $CG_n=G_n(g, zg)$, where $g=\frac{1-\sqrt{1-4z}}{2z}$.
Suppose that $C_1(z)$ and $C_2(z)$ are the generating functions for
the first and the second columns of the binary Riordan matrix ${\cal
B}(zg,zg)_n$. We have that $C_1(z)=zg=\frac{1-\sqrt{1-4z}}{2}$ and
$C_2(z)=zgf=z^2
g^2=\frac{2-2\sqrt{1-4z}-4z}{4}=\frac{1-\sqrt{1-4z}}{2}-z$. From
this, using the fact that $12\in E(CG_n)$, we obtain
$N_{CG_n}(1)\backslash \{2\}=N_{CG_n}(2)\backslash \{1\}$. Now, by
Lemma~\ref{lem5}, the desired result is obtained.
\end{pf}

The following result is an immediate corollary of  Lemmas~\ref{lem5}
and~\ref{dg7}.

\begin{theorem}
 For a Pascal graph $PG_n$, we have
 \begin{itemize}
 \item[{\rm(i)}] If $n=2^p+1$, then  $-1$ is an eigenvalue of $PG_n$ with multiplicity at least $2$ and with eigenvectors
 $X=[\underbrace{1,0,\ldots,0}_{\frac{n-1}{2}},-1,0\ldots,0]^T$ and
 $Y=[1,0,\ldots,0,-1]^T$.
 \item[{\rm(ii)}] If $n\neq 2^p+1$, then  $-1$ is an eigenvalue of $PG_n$ with eigenvector
 $X=[\underbrace{1,0,\ldots,0}_{2^{p}},-1,0,\ldots,0]^T$.
 \end{itemize}
\end{theorem}

\begin{lemma}[\cite{Das3}]\label{lem7}
Let $G$ be a graph of order $n$. If $uv\notin E(G)$ and
$N_G(u)=N_G(v)$ then $0$ is an eigenvalue of $G$ with eigenvector
$$X=[0,\ldots,0,1\,,0, \ldots,0,-1\,,0,\ldots,0]^T,$$ where the
entries $1$ and $-1$ correspond to the vertices $u$ and $v$,
respectively.
\end{lemma}

\begin{theorem}\label{thm10}
Let $2^p\le n\le 2^p+2$. Then $0$ is an eigenvalue of a Pascal graph
$PG_n$ with eigenvector $X$, where the entries of $X$ in positions
$2$ and $2^p$  are $1$ and $-1$, respectively, and the other entries
are $0$s.
\end{theorem}
\begin{pf}
From \cite{CKM}, we have $N_{PG_n}(2)=V_o$  and
$N_{PG_n}(2^p)=\{j\in V_o\,|\, 1\le j \le 2^p+1 \}$. Since $2^p\le
n\le 2^p+2$, we have $N_{PG_n}(2)=N_{PG_n}(2^p)$. Finally, by the
fact that $2(2^p)\notin E(PG_n)$ and Lemma \ref{lem7}, we obtain the
desired.\end{pf}

Now, we are ready to present some results on integral Laplacian
eigenvalues of Pascal and Catalan graphs. The following result is obtained by Lemma 3.1 in \cite{Das}.
\begin{lemma} \label{lem52}
Let $G$ be a graph of order $n$. If $uv\notin E(G)$ and
$N_G(u)=N_G(v)$, then $d_G(u)=|N_G(u)|$ is a Laplacian eigenvalue of
$G$ with eigenvector
$$X=[0,\ldots,0,1\,,0, \ldots,0,-1\,,0,\ldots,0]^T,$$ where the entries
$1$ and $-1$ correspond to the vertices $u$ and $v$, respectively.
\end{lemma}

\begin{theorem}
Let $2^p\le n\le 2^p+2$. Then $\left\lceil\frac{n}{2}\right\rceil$
is a Laplacian eigenvalue of a Pascal graph $PG_n$  with eigenvector
$X$, where the entries in positions $2$ and $2^p$  are $1$ and $-1$,
respectively, and the other entries are $0$s.
\end{theorem}

\begin{pf}
By the proof of Theorem \ref{thm10}, we have
$N_{PG_n}(2)=N_{PG_n}(2^p)$. Since $2(2^p)\notin E(PG_n)$, the required result is obtained by Lemma~\ref{lem52}.
\end{pf}

The following result is obtained by Theorem 3.3 in \cite{Das2}.
\begin{lemma} \label{lem33}
Let $G$ be a graph of order $n$. If $uv\in E(G)$ with
$N_G(u)\backslash  \{v\}=N_G(v)\backslash  \{u\}$, then $d_G(u)+1$
is a Laplacian eigenvalue of $G$ with eigenvector
$$X=[0,\ldots,0,1\,,0, \ldots,0,-1\,,0,\ldots,0]^T,$$ where the entries
$1$ and $-1$ correspond to the vertices $u$ and $v$, respectively.
\end{lemma}

The following result is related to Laplacian eigenvalues of Catalan
graphs.
\begin{theorem}
Let $n\ge 2$. Then $\lceil\log_2 n\rceil+1$ is a Laplacian
eigenvalue of a Catalan graph $CG_n$ with eigenvector
$X=[1,-1,0,\ldots,0]^T$.
\end{theorem}
\begin{pf}
By the arguments used in the proof of Theorem \ref{thm13}, we have
$$N_{CG_n}(1)\backslash \{2\}=N_{CG_n}(2)\backslash \{1\}.$$
Now, by Lemmas~\ref{lem4} and \ref{lem33}, the desired result is
obtained.
\end{pf}

The following result follows directly from Lemma~\ref{dg7} and
Theorem 3.3 in~\cite{Das2}.
\begin{theorem}\label{thm34}
 For a Pascal graph $PG_n$, we have
 \begin{itemize}
 \item[\rm(i)] If $n=2^p+1$, then  $n$ is a Laplacian eigenvalue of $PG_n$ with multiplicity at least $2$ and with eigenvectors
 $X=[\underbrace{1,0,\ldots,0}_{\frac{n-1}{2}},-1,0\ldots,0]^T$ and
 $Y=[1,0,\ldots,0,-1]^T$.\\
 \item[\rm(ii)] If $n\neq 2^p+1$, then  $n$ is a Laplacian eigenvalue of $PG_n$ with eigenvector
 $$X=[\underbrace{1,0,\ldots,0}_{2^{p}},-1,0,\ldots,0]^T.$$
 \end{itemize}
\end{theorem}

We note that Theorem 4.3 in \cite{CKM} gives a better result than
Theorem~\ref{thm34} about the number of Laplacian eigenvalues equal
to $n$, but without any information about the eigenvectors.

\section{The inertias of Riordan graphs}\label{inertia-sec}

 In this section, we present some results on the positive inertia and the negative inertia of Riordan graphs. We begin with two lemmas.

 \begin{lemma} \label{thm2}
 Let $G_n=G_n(g,f)$ be a Riordan graph with $n$ vertices. Then, the following are equivalent:
 \begin{itemize}
 \item[{\rm(i)}]  $[z^{2k}] (gf)\equiv 0$ for $1\le k\le \left\lfloor\frac{n-2}{2}\right\rfloor$.
\item[{\rm(ii)}]  $2j\notin E(G_n)$ for $j\in\{4, 6, \ldots,2\left\lfloor\frac{n}{2}\right\rfloor\}$.
\item[{\rm(iii)}]  $G_n$ is o-decomposable. 
\end{itemize}
\end{lemma}

\begin{pf}
 The statements (i) and (ii) are clearly equivalent. We have (i) holds if and only if  $$\left(\frac{gf}{z}\right)^{\prime}(\sqrt{z})=\sum\limits_{k=1}^{\lfloor\frac{n-2}{2}\rfloor}
 \left(\sum_{i=1}^{2k}g_{2k-i} f_i\right)\, z^{2k}=\sum\limits_{k=1}^{\lfloor\frac{n-2}{2}\rfloor} \left([z^{2k}] (gf)\right)\, z^{2k} \equiv 0.$$
 By Theorem~\ref{e:th}, this is true if and only if  $Y=O_{\left\lceil\frac{n}{2} \right\rceil}$, i.e.  $\langle V_e\rangle\cong N_{\left\lfloor \frac{n}{2}\right\rfloor}$, and $G_n$ is e-decomposable by definition.
 \end{pf}

 \begin{lemma} \label{thm1}
 Let $G_n= G_n(g,f)$ be a Riordan graph with $n$ vertices. Then, the following are equivalent:
 \begin{itemize}
 \item[{\rm(i)}] $[z^{2i-1}] g\equiv 0$ for $1\le i\le \lfloor\frac{n-1}{2}\rfloor$.
 \item[{\rm(ii)}] $1j\notin E(G_n)$ for $j\in\{3, 5, \ldots,2\lceil\frac{n}{2}\rceil-1\}$.
\item[{\rm(iii)}] $G_n$ is e-decomposable. 
 \end{itemize}
 \end{lemma}

\begin{pf}
 The statements (i) and (ii) are clearly equivalent. Further, we have that (i) holds if and only if  $$g^{\prime}(\sqrt{z})=\sum\limits_{i=1}^{\lfloor\frac{n-1}{2}\rfloor} g_{2i-1} z^{2i-1}\equiv 0.$$
 By Theorem~\ref{e:th}, this is true if and only if  $X=O_{\left\lceil\frac{n}{2} \right\rceil}$, i.e.  $\langle V_o\rangle\cong N_{\left\lceil \frac{n}{2}\right\rceil}$, and $G_n$ is e-decomposable by definition.
 \end{pf}

  The following result is obtained by Lemmas \ref{thm1} and \ref{thm2}.

\begin{corollary}\label{cor1}
 Let $G_n$ be a Riordan graph with $n$ vertices.
 \begin{itemize}
 \item[{\rm(i)}] If $G_n=G_{n}(g,zg)$ is of the Bell type, then $G_n$ is o-decomposable.  
 \item[{\rm(ii)}]  If $G_n=G_{n}(f',f)$ is of the derivative type, then $G_n$ is e-decomposable.
\end{itemize}
\end{corollary}

 Now, using the lemmas above, we obtain the following results on the positive and negative inertias of Riordan graphs.
\begin{theorem} \label{th1}
 Let  $G_n$ be an o-decomposable Riordan graph.
 Then,  $$\max\{n^+,\,n^-\}\le \left\lceil \frac{n}{2} \right\rceil.$$
\end{theorem}
\begin{pf} By Lemma \ref{thm2}, the null graph $N_{\left\lfloor\frac{n}{2} \right\rfloor}$ is an induced subgraph of $G_{n}$. This fact, and the interlacing property in Lemma~\ref{lem31}, gives
\begin{equation}\label{eq1}
\lambda_i(G_n)\ge \lambda_i\left(N_{\left\lfloor\frac{n}{2} \right\rfloor}\right)\ge
\lambda_{\left\lceil \frac{n}{2} \right\rceil+i}(G_n)~~~~~\mbox{for}~ 1\le
i\le \left\lfloor\frac{n}{2}\right\rfloor.
\end{equation}  Putting $i=\left\lfloor\frac{n}{2}\right\rfloor$ in the first inequality of (\ref{eq1}), we obtain $$\lambda_{\left\lfloor\frac{n}{2}\right\rfloor}(G_n)\ge \lambda_{\left\lfloor\frac{n}{2}\right\rfloor}\left(N_{\left\lfloor\frac{n}{2} \right\rfloor}\right)=0,$$ that is, $n^+ +\eta\ge \left\lfloor\frac{n}{2}\right\rfloor$, which implies that $n^-\le \left\lceil \frac{n}{2} \right\rceil.$ Now, putting $i=1$ in the second inequality of (\ref{eq1}), we obtain
$$\lambda_{\left\lceil \frac{n}{2} \right\rceil+1}(G_n)\le \lambda_1\left(N_{\left\lfloor\frac{n}{2} \right\rfloor}\right)=0,$$ that is, $n^-+\eta\ge \left\lfloor\frac{n}{2} \right\rfloor$, which implies $n^+ \le \left\lceil \frac{n}{2} \right\rceil$. We are done.
\end{pf}

\begin{corollary} Let  $G_n$ be an o-decomposable Riordan
graph. If $n$ is even and $\det(G_n)\neq 0$, then
$n^+=n^-=\frac{n}{2}$.
\end{corollary}

\begin{corollary} \label{cor3}
 Let $G_n=G_n(g, zg)$ be a Riordan graph of the Bell type with $n$ vertices. Then,
 \begin{equation} \label{equ2}
 \max\{n^+,\,n^-\}\le \left\lceil \frac{n}{2} \right\rceil.
 \end{equation} Moreover, if $n$ is even and $\det(G_n)\neq 0$, then $n^+=n^-=\frac{n}{2}$.
\end{corollary}

The following result is obtained using similar arguments as those in the proofs
of Theorem~\ref{th1} and~Lemma \ref{thm1}.

\begin{theorem} \label{th2-1}
 Let  $G_n$ be an e-decomposable Riordan graph. Then,  $$\max\{n^+,\,n^-\}\le \left\lfloor \frac{n}{2} \right\rfloor.$$
\end{theorem}

Next, we state two corollaries of Theorem~\ref{th2-1}.

\begin{corollary}
  Let  $G_n$ be an e-decomposable Riordan graph. We have the following:
 \begin{itemize}
 \item[{\rm(i)}]   If $n$ is odd, then $\det(G_n)=0$.
 \item[{\rm(ii)}]  If $\det(G_n)\neq 0$, then $n$ is even and $n^+=n^-=\frac{n}{2}$.
 \end{itemize}
\end{corollary}

\begin{corollary}\label{cor4}
 Let $G_n=G_{n}(f',f)$ be a Riordan graph of the derivative type. We have the following:
 \begin{itemize}
 \item[{\rm(i)}] $\max\{n^+,\,n^-\}\le \lfloor \frac{n}{2} \rfloor.$
 \item[{\rm(ii)}]   If $n$ is odd, then $\det(G_n)=0$.
 \item[{\rm(iii)}]  If $\det(G_n)\neq 0$, then $n$ is even and $n^+=n^-=\frac{n}{2}$.
 \end{itemize}
\end{corollary}

 \begin{example} In what follows, we used Sage \cite{sage} to calculate $n^+$ and $n^-$.
 \begin{itemize}
 \item[{\rm(i)}]  The Pascal graph $PG_{10}$ is a Riordan graph of the Bell type. We have $n^+(PG_{10})+1= n^-(PG_{10})=5$. This is consistent with Corollary~\ref{cor3} and gives the equality in (\ref{equ2}).
 \item[{\rm(ii)}]   The graph $G_{10}\left(\frac{1}{1+z^2}, \frac{z}{1+z}\right)$ is a Riordan graph of the derivative type. We have\\ $n^+\left(G_{10}\left(\frac{1}{1+z^2}, \frac{z}{1+z}\right)\right)= n^-\left(G_{10}\left(\frac{1}{1+z^2}, \frac{z}{1+z}\right)\right)=5$. This is consistent with Corollary~\ref{cor4} and gives the equality in (i).
 \item[{\rm(iii)}]  $n^+\left(G_{16}\left(1+z^3, \frac{z}{1+z}\right)\right)= 6,~n^-\left(G_{16}(1+z^3, \frac{z}{1+z})\right)=10$.
 \item[{\rm(iv)}]   $n^+(K_n)=1,~n^-(K_n)=n-1$.
  \end{itemize}
 Clearly, $G_{16}\left(1+z^3, \frac{z}{1+z}\right)$ and $K_n$ do not satisfy the conditions in Lemmas~\ref{thm1} and~\ref{thm2}.
 \end{example}

  Using Sage \cite{sage}, we observe that $n^-(PG_i)=\lceil \frac{i}{2}\rceil$ for $2\le i\le 200$. Moreover, by Lemma~\ref{check}, any Riordan graph $G_n=G_n(g,zg)$ with an even function $g$ is bipartite and thus $n^+(G_n)=n^-(G_n)$. Based on these observation, we state the following problem.

\begin{problem}
 Is it true that $n^+\le n^-$ for any Riordan graph $G_n(g,zg)$?
\end{problem}

Let $H_n$ be the bipartite graph obtained from a Riordan graph $G_n(g,f)$ by deleting all edges in
$\langle V_o\rangle$ and $\langle V_e\rangle$.
 The following result is obtained directly from Lemma~\ref{lem41}.

\begin{theorem} Let $G_n=G_n(g ,f)$ be an o-decomposable, or an e-decomposable, Riordan graph. Then,
$$\min{\left(n^+(G_n), n^-(G_n)\right)}\ge \frac{1}{2} {\rm rank} (H_n)=n^+(H_n)=n^-(H_n).$$ In particular, if $n$ is even and the corresponding bipartite graph is nonsingular, then
$n^-(G_n)=n^+(G_n)=n/2$ and $\eta(G_n)=0$.
\end{theorem}

For a graph $G$, we denote the number of eigenvalues of $G$ located
in the interval $I$ by $m_G(I)$. The following results give
bounds on $m_{\overline{G}_n}(-\infty,-1]$.
\begin{theorem} Let $G_n=G_n(g ,f)$ be an o-decomposable Riordan graph and $\overline{G}_n$ be its complement. Then, $$n^-(\overline{G}_n)\ge m_{\overline{G}_n}(-\infty,-1]\ge \left\lfloor\frac{n}{2}\right\rfloor-1.$$
\end{theorem}

\begin{pf}
From the definition of an o-decomposable Riordan graph, we have
$\overline{K}_{\lfloor n/2 \rfloor}$ is an induced subgraph of $G_n$
and, consequently,  the complete graph $K_{\lfloor n/2 \rfloor}$ is
an induced subgraph of $\overline{G}_n$. Now, by the well known
interlacing theorem, we obtain
$$-1=\lambda_2(K_{\lfloor n/2 \rfloor})\ge \lambda_{n-\lfloor n/2
\rfloor+2}(\overline{G}_n)=\lambda_{\lceil
n/2\rceil+2}(\overline{G}_n),$$ which gives the desired result.
\end{pf}

By a similar argument, one can easily obtain the following result.
\begin{theorem} Let $G_n=G_n(g ,f)$ be an e-decomposable Riordan graph and $\overline{G}_n$ be its complement. Then, $$n^-(\overline{G}_n)\ge m_{\overline{G}_n}(-\infty,-1]\ge \left\lceil\frac{n}{2}\right\rceil-1.$$
\end{theorem}

\section{The nullity of Riordan graphs}\label{nullity-sec}

 In this section, we present some results on the nullity of Riordan graphs.

\begin{theorem}\label{adjacency}
Let ${\cal A}_n$ be the adjacency matrix of an o-decomposable
Riordan graph $G_n(g,f)$ of the form $P {\footnotesize
\left[
                     \begin{array}{cc}
                         X_{\left\lceil \frac{n}{2} \right\rceil}& B\\[3mm]
                       B^T& $O$\\
                     \end{array}
                   \right]} P^T$.
Then
 $$\eta(B)\le \eta({\cal A}_n)\le \left\{
     \begin{array}{ll}
        2\eta(B) & \hbox{if n is even} \\[3mm]
        2\eta(B)+1 & \hbox{if n is odd.}
     \end{array}
   \right.
 $$
\end{theorem}
\begin{pf}
We have $${\rm rank}(B)+{\rm rank}(B^T)={\rm rank}\left(\left[
                     \begin{array}{cc}
                         O_{\left\lceil \frac{n}{2} \right\rceil}& B\\
                       B^T& $O$\\
                     \end{array}
                   \right]\right)\le {\rm rank}({\cal A}_n)\le {\rm rank}\left(\left[
                     \begin{array}{c}
                         X_{\left\lceil \frac{n}{2} \right\rceil}\\
                       B^T\\
                     \end{array}
                   \right]\right)+{\rm rank}(B).$$
Since ${\rm rank}(B)={\rm rank}(B^T)$, from the above, we obtain
$$2\,{\rm rank}(B)\le {\rm rank}({\cal A}_n)\le \left\lceil\frac{n}{2}\right\rceil+{\rm rank}(B).$$ We can now apply the rank-nullity theorem to obtain
$$2\,\left(\left\lfloor\frac{n}{2}\right\rfloor-\eta(B)\right)\le n-\eta({\cal A}_n)\le \left\lceil\frac{n}{2}\right\rceil+\left(\left\lfloor\frac{n}{2}\right\rfloor-\eta(B)\right)~\Rightarrow~\eta(B)\le \eta({\cal A}_n)\le 2\eta(B)+n-2\left\lfloor\frac{n}{2}\right\rfloor.$$ This completes our proof.
\end{pf}

\begin{corollary}
Let $G_n(g,f)$ be an o-decomposable Riordan graph of even order $n$
with the adjacency matrix ${\cal A}_n$ of the form in Theorem
\ref{adjacency}. Then  $\det{B}\neq 0$ if and only if $\det{{\cal
A}_n}\neq 0$.
\end{corollary}

\begin{corollary}
Let $G_n(g,f)$ be an o-decomposable Riordan graph of odd order $n$
with the adjacency matrix ${\cal A}_n$ of the form in Theorem
\ref{adjacency}. If  $\eta{(B)}=0$ then $\eta{({\cal A}_n)}=0$ or
$1$.
\end{corollary}

 When the graph $G_n(g, zg)$ is proper (recall Definition~\ref{def-prop-graph}), we have the following result on the matrix $P^{T}{\cal A}_n P$, where ${\cal A}_n={\cal A}_n(g,zg)$.
This result is obtained by Lemma~\ref{lem12}.

 \begin{corollary} \label{rm1}
  Let $[z^0]g\equiv 1$ and ${\cal A}_n={\cal A}_n(g, zg)$. Then in the matrix $P^{T}{\cal A}_n P$ in Lemma~\ref{lm3}, we have the following.
  \begin{itemize}
\item[{\rm(i)}]  If $n$ is odd, then  $$b_{i,i}=1~~~~\mbox{for}~~1\le i \le \frac{n-1}{2},~~~\mbox{and}~~~~~r_{j,(j-1)}=1~~~~\mbox{for}~~2\le j \le \frac{n+1}{2}.$$
\item[{\rm(ii)}]  If $n$ is even, then $$b_{i,i}=1~~~~\mbox{for}~~1\le i \le \frac{n}{2},~~~\mbox{and}~~~~~r_{j,(j-1)}=1~~~~\mbox{for}~~2\le j \le \frac{n}{2}.$$
 \end{itemize}
 \end{corollary}

 \begin{theorem}\label{lem3}
 Let $[z^0]g\equiv 1$ and let
 $$
 {\cal A}_n:= P {\footnotesize \left[
                     \begin{array}{cc}
                         {\cal A}_{\left\lceil \frac{n}{2} \right\rceil}& B\\[3mm]
                       B^T& $O$\\
                     \end{array}
                   \right]} P^T
                   $$
                    be the adjacency matrix of an io-decomposable Riordan graph $G_n(g,zg)$ of the Bell type. For any $n\ge 2$,  we have $${\rm rank}\left(\left[
                     \begin{array}{c}
                         {\cal A}_{\left\lceil \frac{n}{2} \right\rceil}\\[3mm]
                       B^T\\
                     \end{array}
                   \right]\right)=\left\lceil \frac{n}{2}\right\rceil. $$
 \end{theorem}
 \begin{pf}
 We consider the cases of even and odd $n$ separately. \\
(i) Let $n=2k$ and $$\left[
                    \begin{array}{c}
                      {\cal A}_k \\[2mm]
                      B^T \\
                    \end{array}
                  \right]=\left[
                            \begin{array}{cccccc}
                              {\bf r}_1&\cdots&{\bf r}_k&{\bf b}_1 &\cdots &{\bf b}_k \end{array} \right]^T
$$
where ${\bf r}_i$ is the $i$th row of ${\cal A}_k$ and ${\bf b}_i$  is  the $i$th row of $B^T$.
We consider the $k \times k$ matrix
$$X_k=\left[\begin{array}{cccc}{\bf b}_1-{\bf r}_1 &\cdots & {\bf b}_{k-1}-{\bf r}_{k-1}
&{\bf b}_k \end{array}\right]^T.
$$ It follows from Lemma \ref{lm3} and Remark \ref{rm1} that  $X_k$ is a unit lower triangular matrix. Thus, $$ {\rm rank}\left(\left[
                     \begin{array}{c}
                        {\cal A}_{\left\lceil \frac{n}{2} \right\rceil}\\[3mm]
                       B^T\\
                     \end{array}
                   \right]\right)= {\rm rank}(X_k)=k=\left\lceil \frac{n}{2}\right\rceil.$$
  \noindent
(ii) Let $n=2k+1$. Note that $ {\cal A}_{\left\lceil \frac{n}{2}
\right\rceil}= {\cal A}_{k+1}$ and $B^T$ is the $k\times (k+1)$
matrix. Using the same notation as above, we have
 $$
 \left[
                    \begin{array}{c}
                      {\cal A}_{k+1} \\[2mm]
                      B^T \\
                    \end{array}
                  \right]=\left[
                            \begin{array}{cccccc}
                              {\bf r}_1&\cdots&{\bf r}_{k+1}&{\bf b}_1 &\cdots &{\bf b}_k \end{array} \right]^T.
$$
We consider the $(k+1) \times (k+1) $ matrix
$$
Y_{k+1}=\left[\begin{array}{cccc}{\bf b}_1-{\bf r}_1 &\cdots & {\bf
b}_k-{\bf r}_k &{\bf r}_k \end{array}\right]^T.
$$
It follows from Lemma \ref{lm3} and Remark \ref{rm1} that $Y_{k+1}$
is also a unit lower triangular matrix.
                                     Thus, $$ {\rm rank}\left(\left[
                     \begin{array}{c}
                        {\cal A}_{\left\lceil \frac{n}{2} \right\rceil}\\[3mm]
                       B^T\\
                     \end{array}
                   \right]\right)= {\rm rank}(Y_{k+1})=k+1=\left\lceil \frac{n}{2}\right\rceil,$$ which completes our proof.
 \end{pf}

 Let $X=[x_1, x_2, \ldots, x_n]^T$ be an eigenvector of the matrix ${\cal A}_n$ with eigenvalue $0$. We define $X_o$ and $X_e$ to be the vectors of all eigencomponents $x_i$, respectively, with odd and even indices $i$,  that is,
  $$X_o=\left[x_1, x_3, \ldots, x_{2\left\lceil\frac{n}{2}\right\rceil-1}\right]^T~~~\mbox{and}~~X_e=\left[x_2, x_4, \ldots, x_{2\left\lfloor\frac{n}{2}\right\rfloor}\right]^T.$$

 \begin{lemma} \label{th2}
Let ${\cal A}_n$ be the adjacency matrix of a singular
io-decomposable Riordan graph $G_n(g,zg)$ of the form in Theorem
\ref{lem3}. Then
$$X_o=O_{\lceil \frac{n}{2}\rceil}~~~\mbox{if and only
if}~~~BX_e=O_{\left\lceil \frac{n}{2}\right\rceil}.$$
\end{lemma}
\begin{pf} By Theorem \ref{lem12} and our assumptions, we have
 \begin{eqnarray}
 A_n\, X=\left[
 \begin{array}{cc}
 {\cal A}_{\left\lceil \frac{n}{2}\right\rceil} & B \\[3mm]
 B^T & O
 \end{array}
 \right]\, \left[
 \begin{array}{c}
 X_o \\[3mm]
 X_e
 \end{array}
 \right]=O_n,
 \end{eqnarray}
 \begin{eqnarray}
 {\cal A}_{\left\lceil \frac{n}{2}\right\rceil}\,X_o+BX_e &=& O_{\left\lceil \frac{n}{2}\right\rceil},\nonumber\\
 B^T\, X_o  &=& O_{\lfloor \frac{n}{2}\rfloor}.\nonumber
 \end{eqnarray}
If $X_o=O_{\left\lceil \frac{n}{2}\right\rceil}$, obviously, we obtain
$BX_e=O_{\left\lceil \frac{n}{2}\right\rceil}$. Now, we assume that
$BX_e=O_{\left\lceil \frac{n}{2}\right\rceil}$. Then, $D\,X_o=O_n,$ where
$$D=\left[
 \begin{array}{c}
 {\cal A}_{\left\lceil \frac{n}{2}\right\rceil} \\[3mm]
 B^T
 \end{array}
 \right].$$ Since, by Theorem \ref{lem3}, the matrix $D$ is full column rank, we obtain $X_o=O_{\left\lceil \frac{n}{2}\right\rceil}$, as desired.
\end{pf}

Using Sage \cite{sage}, we come up with the following conjecture.

\begin{conjecture}\label{conj1}
Let $G_n(g,zg)$ be a singular io-decomposable Riordan graph. Then
$X_o=O_{\left\lceil\frac{n}{2}\right\rceil}$.
\end{conjecture}

In what follows, we show that Conjecture~\ref{conj1} is equivalent
to a condition on $\eta(G_n)$.

\begin{theorem}\label{thm21}
Let ${\cal A}_n$ be the adjacency matrix of a singular
io-decomposable Riordan graph $G_n(g,zg)$ of the form in Theorem
\ref{lem3}. Then  $X_o=O_{\left\lceil\frac{n}{2}\right\rceil}$ if
and only if $\eta(G_n)=\eta(B)$.
\end{theorem}

\begin{pf} First, assume that $X_o=O_{\lceil\frac{n}{2}\rceil}$.  Then, by
Lemma \ref{th2}, this is equivalent to
$BX_e=O_{\left\lceil\frac{n}{2}\right\rceil}$. This leads to
$\eta(G_n)=\eta(B)$.

Now, we assume that $\eta(G_n)=\eta(B)$, that is, that ${\rm rank}(
{\cal A}_n)=\lceil\frac{n}{2}\rceil+{\rm rank}(B)$. Since $X$ is an
eigenvector of $G_n$ with eigenvalue $0$, we have
$$\left[
                     \begin{array}{cc}
                         {\cal A}_{\left\lceil \frac{n}{2} \right\rceil}& B\\[3mm]
                       B^T& O\\
                     \end{array}
                   \right] \left[
                                                                   \begin{array}{c}
                                                                     X_o \\
                                                                     X_e\\
                                                                   \end{array}
                                                                 \right]=O_n.$$
Suppose that ${\rm rank}(B)=r_b$.  Using row operations on
$B$, one can easily obtain $\left(
                                                              \begin{array}{c}
                                                                B_{r_b} \\
                                                                O \\
                                                              \end{array}
                                                            \right)
$, where $B_{r_b}$ is an $r_b \times \lfloor\frac{n}{2}\rfloor$
matrix containing  $r_b$ linearly independent rows of $B$. Using the
same row operations on $P^T {\cal A} P$, we have
$$\left[
    \begin{array}{c|c}
      A_1 & B_{r_b}\\
\hline
      A_2& O\\
\hline
      B^T& O \\
    \end{array}
  \right]\left[
                                                                   \begin{array}{c}
                                                                     X_o \\
                                                                     X_e\\
                                                                   \end{array}
                                                                 \right]=O_n,$$
where $A_2$ is an $\left(\left\lceil\frac{n}{2}\right\rceil-r_b\right)\times\left\lceil\frac{n}{2}\right\rceil$ matrix. From this, we obtain $\left[
                           \begin{array}{c}
                             A_2 \\
                             B^T \\
                           \end{array}
                         \right] X_o=O_{k}$, where $k=\lfloor\frac{n}{2}\rfloor+\lceil\frac{n}{2}\rceil-r_b$. To prove $X_o=O_{\left\lceil\frac{n}{2}\right\rceil}$, it suffices to show
that \begin{equation}\label{moj4} {\rm
rank}\left(\left[\begin{array}{c}
                             A_2 \\
                             B^T \\
                           \end{array}
                         \right] \right)=\left\lceil\frac{n}{2}\right\rceil.
\end{equation}

We have $$\left\lceil\frac{n}{2}\right\rceil+{\rm rank}(B)={\rm rank}\left[
    \begin{array}{c|c}
      A_1 & B_{r_b}\\
\hline
      A_2& O\\
\hline
      B^T& O \\
    \end{array}
  \right]\le {\rm rank}(B_r)+{\rm rank}\left(\left[
                           \begin{array}{c}
                             A_2 \\
                             B^T \\
                           \end{array}
                         \right]\right),$$ which gives the desired result
in (\ref{moj4}). Hence, the proof is done.
\end{pf}

Applying Theorem~\ref{thm21}, we can state a different version of Conjecture~\ref{conj1}.

\begin{conjecture}\label{conj2}
Let ${\cal A}_n$ be the adjacency matrix of a singular
io-decomposable Riordan graph $G_n(g,zg)$ of the form in Theorem
\ref{lem3}. Then, $\eta(G_n)=\eta(B)$.
\end{conjecture}

\begin{lemma}[\cite{Ber}]\label{moj5}
Let $A, B, C$ and $D$ be, respectively, real $n\times n$, $n\times m$,
$l \times n$ and $l \times m$ matrices, and let $A$ be nonsingular.
Then
$${\rm rank}\left[
        \begin{array}{cc}
          A & B \\
          C & D \\
        \end{array}
      \right]=n+{\rm rank}(D-CA^{-1} B).$$
\end{lemma}

\begin{theorem}\label{thm14}
Let ${\cal A}_n$ be the adjacency matrix of an io-decomposable
Riordan graph $G_n(g,zg)$ of the form in Theorem \ref{lem3}. Then
$\eta(G_n)=\eta(B^T M^{-1} B)$, where $M= {\cal
A}_{\frac{n}{2}}-B^T$ for even $n$, and $$ M=\left[
\begin{array}{cccc}{\bf r}_1-{\bf b}_1 &\cdots &{\bf r}_{k-1}- {\bf b}_{k-1}&{\bf r}_k-{\bf b}_{k-1}\end{array}\right]^T$$ for odd $n$,
$k=\frac{n+1}{2}$, where ${\bf r}_i$ is the $i$th row of
${\cal A}_k$ and ${\bf b}_i$  is  the $i$th row of $B^T$.
\end{theorem}
\begin{pf}
We know that the rank of any matrix is fixed under row and column
operations. Thus, we have
$${\rm rank}({\cal A}_n)={\rm rank}\left[
                     \begin{array}{cc}
                         {\cal A}_{\left\lceil \frac{n}{2} \right\rceil}& B\\[3mm]
                       B^T& $O$\\
                     \end{array}
                   \right]={\rm rank}\left[
                     \begin{array}{cc}
                         M_{\left\lceil \frac{n}{2} \right\rceil}& B\\[3mm]
                       B^T& $O$\\
                     \end{array}
                   \right].$$ By Lemma \ref{lm3}, the square matrix $M$
is a lower triangular matrix with all entries on the main diagonal
equal to $-1$. Therefore, $M$ is invertible, and the required result
is given by Lemma~\ref{moj5}.
\end{pf}

 By Theorem \ref{thm14}, it is worth noticing that Conjecture
\ref{conj2} is true for any io-decomposable Riordan graph
$G_n(g,zg)$ if and only if $\eta(B)=\eta(B^T M^{-1} B)$.

 \section{Determinants of Riordan graphs}\label{det-Catalan-sec}

 In this section, we study determinants of  o-decomposable and e-decomposable Riordan graphs, as well
as determinants of Catalan graphs.

We begin with Schur's formula for determinant
of a block matrix.
\begin{lemma}[\cite{Ber}]\label{schur}
Let $A, B, C$ and $D$ be real $n\times n$ matrices. Then
$${\rm det}\left[
        \begin{array}{cc}
          A & B \\
          C & D \\
        \end{array}
      \right]=\left\{
                \begin{array}{ll}
                  {\rm det}(DA-CB), & \mbox{if } \hbox{AB=BA;} \\[2mm]
                  {\rm det}(AD-BC), & \mbox{if } \hbox{DC=CD.}
                \end{array}
              \right.
$$
\end{lemma}

Recall that $H_n$ is the bipartite graph obtained from a Riordan graph $G_n(g,f)$ by deleting all edges in
$\langle V_o\rangle$ and $\langle V_e\rangle$.

\begin{theorem}\label{det}
Let $G_n=G_n(g,f)$ be an o-decomposable, or e-decomposable, Riordan
graph of even order $n$. Then, $${\rm det}(G_n)=({\rm det}(B))^2,$$
where $B$ is given in (\ref{e:bm}).
\end{theorem}
\begin{pf}
First, suppose that  ${\cal A}_n$ is the adjacency matrix of an
o-decomposable Riordan graph $G_n$. We have
$${\cal A}_n:=P {\footnotesize \left[
                     \begin{array}{cc}
                         X_{\left\lceil \frac{n}{2} \right\rceil}& B\\[3mm]
                       B^T& $O$\\
                     \end{array}
                   \right]} P^T.$$
Since $n$ is even and $B$ is a square matrix, by Lemma~\ref{schur}, we obtain $${\rm det}(G_n)={\rm det}({\cal A}_n)={\rm
det}\left[
                     \begin{array}{cc}
                         X_{\left\lceil \frac{n}{2} \right\rceil}& B\\[3mm]
                       B^T& $O$\\
                     \end{array}
                   \right]={\rm det}(BB^T)=({\rm det}(B))^2.$$
This completes the proof for o-decomposable Riordan graphs.

Next, suppose that  ${\cal A}_n$ is the adjacency matrix of an
e-decomposable Riordan graph $G_n$. We have
$${\cal A}_n:=P {\footnotesize \left[
                     \begin{array}{cc}
                         $O$& B\\[3mm]
                       B^T& Y_{\left\lfloor \frac{n}{2} \right\rfloor}\\
                     \end{array}
                   \right]} P^T.$$
Since $n$ is even, $B$ is a square matrix, and thus, by Lemma~\ref{schur}, we obtain $${\rm det}(G_n)={\rm det}({\cal A}_n)={\rm
det}\left[
                     \begin{array}{cc}
                         $O$& B\\[3mm]
                       B^T& Y_{\left\lfloor \frac{n}{2} \right\rfloor}\\
                     \end{array}
                   \right]={\rm det}(B^T B)=({\rm det}(B))^2.$$
This completes the proof for e-decomposable Riordan graphs.
\end{pf}

The following results are obtained by Theorem \ref{det}.

\begin{corollary} \label{det1}
Let $G_n=G_n(g,f)$ be an o-decomposable Riordan graph of even order
$n$. Then ${\rm det}(B)=0$ if and only if ${\rm det}(G_n)=0$, that
is, $G_n$ is a singular graph if and only if $H_n$ is such a graph.  In
particular, if $N_{H_n}(i)=N_{H_n}(j)$ for $i,j\in V_o$, $i\neq j$,
then $\det(G_n)=0$.
\end{corollary}
\begin{pf}
Our first claim is a direct corollary of Theorem~\ref{det}. Now,
suppose that ${\bf b}_i$ and ${\bf b}_j$ are the $i$th and $j$th rows
of the submatrix $B$. From our assumption, we obtain ${\bf b}_i={\bf b}_j$,
where $i,j\in V_o\,(i\neq j)$. This gives ${\rm det}(B)=0$, and hence,
by Theorem~\ref{det}, the proof is complete.
\end{pf}

\begin{corollary} \label{det11}
Let $G_n=G_n(g,f)$ be an o-decomposable Riordan graph of even order
$n$. If $G_n$ has at least two universal vertices, then ${\rm
det}(G_n)=0$.
\end{corollary}

\begin{corollary} \label{det2}
Let $G_n=G_n(g,f)$ be an e-decomposable Riordan graph of even order
$n$. Then, ${\rm det}(B)=0$ if and only if ${\rm det}(G_n)=0$, that
is, $G_n$ is a singular graph if and only if $H_n$ is such a graph. In
particular, if $N_{H_n}(i)=N_{H_n}(j)$ for $i,j\in V_e$, $i\neq j$,
then ${\rm det}(G_n)=0$.
\end{corollary}

\begin{corollary} \label{det21}
Let $G_n=G_n(g,f)$ be an e-decomposable Riordan graph of even order
$n$. If $G_n$ has at least two universal vertices, then ${\rm
det}(G_n)=0$.
\end{corollary}

The following gives an example for Corollary~\ref{det2}.
\begin{example}
Let $G_8=G_8\left(\frac{1}{1-z^2}, \frac{z}{1-z}\right)$. Since $g$ is an even
function, $G_8$ is an e-decomposable Riordan graph. The matrix $B$
for this graph is
$$B=\left[
      \begin{array}{cccc}
        1 & 1 & 1 & 1 \\
        1 & 1 & 0 & 1 \\
        0 & 1 & 1 & 1 \\
        1 & 1 & 1 & 1 \\
      \end{array}
    \right].
$$ Since ${\bf b}_1={\bf b}_4$, ${\rm det}(B)=0$, and then, by Corollary~\ref{det2},  ${\rm det}(G_8)=0$.
\end{example}

We end this section by considering determinants of Catalan graphs
$CG_n$. Let $CM_n$ be the adjacency matrix of $CG_n$. Since $CG_n$
is a Bell type io-decomposable Riordan graph, from Theorem~12 in
\cite{CJKM}, there exists a permutation matrix $P$ such that
 \begin{eqnarray}\label{ctg}
 A_n=P\, (CM_n)\,P^{T}=\left[
 \begin{array}{cc}
 CM_{\left\lceil \frac{n}{2}\right\rceil} & B \\[3mm]
 B^T & O
 \end{array}
 \right].
 \end{eqnarray}
 By the binomial formula for the Catalan matrix $\left({1-\sqrt{1-4z}\over 2z},{1-\sqrt{1-4z}\over 2}\right)$, we have that
\begin{equation} \label{bi1}
 CM_n(i,j)\equiv C(i-2,j-1)=\frac{j}{i-1}\binom{2i-j-3}{i-j-1},~~~~i>j\ge 1.
\end{equation}

 Let ${\bf b}_j$ be the $j$th row of the submatrix $B$.
\begin{lemma} \label{lem2}
 For the block $B$ in the matrix $A_n$, we have
 \begin{equation}
 {\bf b}_2\equiv{\bf b}_3. \label{eq2}
 \end{equation}
 \end{lemma}

\begin{pf}To prove (\ref{eq2}), it suffices to prove that \begin{equation}\label{eq3} CM_n(3,j)= CM_n(5,j)\end{equation} for each even $j$,  $2\le j \le n$.
For $j\in\{2,4\}$, we have $CM_n(3,j)=CM_n(5,j)=1$.
Suppose that
$j=2k$ where $3\le k \le 2\left\lfloor \frac{n}{2}\right\rfloor$. By
(\ref{bi1}), we have
 \begin{equation} \label{eq4}
 C(2k-2,2)=\frac{3\,(2k-3)\,(4k-7)}{5\,(k-2)\,(2k-5)}C(2k-2,4)\equiv \frac{C(2k-2,4)}{k-2}.
 \end{equation}
If $k$ is odd, then by (\ref{bi1}) and (\ref{eq4}),
 $$CM_n(2k,3)\equiv C(2k-2,2)\equiv \frac{C(2k-2,4)}{k-2} \equiv  CM_n(2k,5),$$
which gives us  (\ref{eq3}). Otherwise $k$ is even. With $C_n=\frac{1}{n+1}{2n\choose n}$ denoting the $n$th Catalan number, we have
 \begin{eqnarray}
 C(2k-2,4)&=& \frac{5}{2k-1} \binom{4k-8}{2k-6} \nonumber\\[2mm]
          &=& \frac{5}{s+3} \binom{2s}{s-2} ~~~~\mbox{(letting}\,s:=2k-4)\nonumber\\
          &=&\frac{5\,(s-1)\,s}{(s+1)\, (s+2)\, (s+3)} \binom{2s}{s}\nonumber\\
          &=&\frac{5\,(s-1)\,s}{(s+2)\, (s+3)}\, C_s,  \nonumber\\
         &\equiv& \frac{s}{s+2}\, C_s= \frac{2k-4}{2k-2}\, C_{2k-4} \equiv  (k-2)\, C_{2k-4}\equiv 0.\label{eq5}
 \end{eqnarray}  Now, by (\ref{eq4}), (\ref{eq5}) and Lemma~\ref{lem4}, we obtain
  $$C(2k-2,2)\equiv \frac{C(2k-2,4)}{k-2}\equiv C_{2k-4}\equiv 0$$ completing our proof.
\end{pf}

\begin{theorem}\label{det-cat-thm}
$\det(CG_{2n})=0$ for any $n\ge 3$.
\end{theorem}

\begin{pf} Since the order of the matrix $CM_{2n}$ is $2n$, the block $B$ in (\ref{ctg}) is a square matrix of order $n$. On the other hand, by Lemma \ref{lem2} we have
 ${\bf b}_2\equiv{\bf b}_3$. Thus, ${\rm rank}(B)\le n-1$ and $\det(B)=0$. Now, we obtain the desired result from the fact that $\det(CM_{2n})=(-1)^{2n}$.
\end{pf}

\begin{remark}
 Let $\Gamma=\{11, 13, 15, 23, 33, 51, 61, 63\}$. By Sage \cite{sage}, we observe that for $n\in \Gamma$,
 $\det(CG_{n})=0$.
\end{remark}

 We conclude our paper with the following conjecture.

\begin{conjecture}
We have that $\det{CG_n}=0$ if and only if $n\ge 6$ is even, or $n\in \Gamma$.
\end{conjecture}

 \end{document}